\documentstyle[11pt]{article}

\newtheorem{eg}{{\sc Example}}
\newtheorem{lemma}{{\sc Lemma}}[section]
\newtheorem{theo}{{\sc Theorem}}

\newtheorem{cor}[theo]{{\sc Corollary}}

\font\bbb=msbm10 scaled\magstep1

\newcommand{\RR}{\mbox{\bbb R}}
\newcommand{\ZZ}{\mbox{\bbb Z}}

\begin{document}

\title{Degree-regular triangulations of the double-torus }

\author{Basudeb Datta and Ashish Kumar Upadhyay}
\date{To appear in `{\bf Forum Mathematicum}'}
\maketitle

\vspace{-5mm}

\hrule

\begin{abstract} A connected combinatorial $2$-manifold is
called degree-regular if each of its vertices have the same
degree. A connected combinatorial $2$-manifold is called weakly
regular if it has a vertex-transitive automorphism group.
Clearly, a weakly regular combinatorial $2$-manifold is
degree-regular and a degree-regular combinatorial $2$-manifold of
Euler characteristic $-2$ must contain $12$ vertices.

In 1982, McMullen {\em et al.} constructed a $12$-vertex
geometrically realized triangulation of the double-torus in
$\RR^3$. As an abstract simplicial complex, this triangulation is
a weakly regular combinatorial $2$-manifold. In 1999, Lutz showed
that there are exactly three weakly regular orientable
combinatorial 2-manifolds of Euler characteristic $-2$. In this
article, we classify all the orientable degree-regular
combinatorial 2-manifolds of Euler characteristic $-2$. There are
exactly six such combinatorial 2-manifolds. This classifies all
the orientable equivelar polyhedral maps of Euler characteristic
$-2$.

\end{abstract}

{\small

{\bf AMS classification\,:} 57Q15, 57M20, 57N05.

{\bf Keywords\,:} Triangulations of $2$-manifolds, regular
simplicial maps, combinatorially  \newline \mbox{} \hspace{25mm}
regular triangulations, degree-regular triangulations.

}

\bigskip

\hrule

\section{Introduction}

Recall that an (abstract) {\em simplicial complex} is a
collection of nonempty finite sets (set of {\em vertices}) such
that every nonempty subset of a member is also a member.  For $i
\geq 0$, the elements of size $i+1$ are called the {\em
$i$-simplices} of the simplicial complex. 1-simplices are also
called the {\em edges} of the simplicial complex. For a
simplicial complex $X$, the maximum of $k$ such that $X$ has a
$k$-simplex is called the {\em dimension} of $X$. The set $V(X)$
of vertices of $X$ is called the {\em vertex-set} of $X$. A
simplicial complex $X$ is called finite if $V(X)$ is finite.

If $X$ and $Y$ are two simplicial complexes, then a (simplicial)
isomorphism from $X$ to $Y$ is a bijection $\varphi \colon V(X)
\to V(Y)$ such that for $\sigma \subseteq V(X)$, $\sigma$ is a
simplex of $X$ if and only if $\varphi(\sigma)$ is a simplex of
$Y$. Two simplicial complexes $X, Y$ are called {\em isomorphic}
(denoted by $X \cong Y$) when such an isomorphism exists. We
identify two simplicial complexes if they are isomorphic. An
isomorphism from a simplicial complex $X$ to itself is called an
{\em automorphism} of $X$. All the automorphisms of $X$ form a
group, which is denoted by ${\rm Aut}(X)$.

A simplicial complex $X$ is usually thought of as a prescription
for constructing a topological space (called the {\em geometric
carrier} of $X$ and is denoted by $|X|$) by pasting together
geometric simplices. Formally, $|X|$ is the subspace of
$[0,1]^{V(X)}$ consisting of the functions $f: V(X) \to [0,1]$
such that the support $\{ v \in V (X) : f (v) \neq 0 \}$ is a
simplex of $X$ and $\sum_{v \in V (X)} f(v) = 1$.  If $\sigma $
is a simplex then $|\sigma  |:= \{ f\in |X| : \sum_{v \in \sigma }
f(v) = 1 \}$ is called the {\em geometric simplex} corresponding
to $\sigma $. We say that a simplicial complex $X$ {\em
triangulates} a topological space $Y$ if $Y$ is homeomorphic to
$|X|$. A simplicial complex $X$ is called {\em connected} if
$|X|$ is connected. A 2-dimensional simplicial complex is called
a {\em combinatorial $2$-manifold} if it triangulates a closed
surface. A combinatorial $2$-manifold $X$ is called {\em
orientable} if $|X|$ is an orientable 2-manifold.

For a simplicial complex $X$, an embedding $i \colon |X|
\hookrightarrow \RR^n$ is called {\em geometric} if the image of
the geometric simplices are convex. A triangulation $K$ of a
surface is called {\em geometrically realizable} in $\RR^3$ if
$|K|$ can be embedded geometrically in $\RR^3$.

If $v$ is a vertex of a simplicial complex $X$, then the number
of edges containing $v$ is called the {\em degree} of $v$ and is
denoted by $\deg_X(v)$ (or $\deg(v)$). A simplicial complex is
called {\em neighbourly} if each pair of vertices form an edge.

A {\em pattern} is an ordered pair $(M,  G)$, where $M$ is a
connected closed surface and $ G$ is a finite graph on $M$ such
that each component of $M \setminus G$ is simply connected. The
closure of each component of $M \setminus  G$ is called a {\em
face} of $(M,  G)$. A pattern is called a {\em polyhedral
pattern} if the boundary of each face is a cycle and the
intersection of any two faces is empty, a vertex or an edge (cf.
\cite{eek}). A {\em polyhedral map} is the collection of the
boundary cycles (together with their vertices and edges) of a
polyhedral pattern (cf. \cite{bw}). So, a combinatorial
2-manifold is a polyhedral map. For a polyhedral map $K$,
consider the polyhedral map $\widetilde{K}$ (called the {\em
dual} of $K$) whose vertices are the faces of $K$ and $(F_1,
\dots, F_m)$ is a face of $\widetilde{K}$ if $F_1, \dots, F_m$
have a common vertex and $F_1 \cap F_2, \dots, F_{m-1} \cap F_m,
F_{m} \cap F_1$ are edges. For a combinatorial definition of
polyhedral maps see \cite{ms}.

A polyhedral map is called {\em $\{p, q\}$-equivelar} if each
face is a $p$-cycle and each vertex is in $q$ faces. A polyhedral
map is called {\em equivelar} if it is $\{p, q\}$-equivelar for
some $p$, $q$ (cf. \cite{ms, msw}). Clearly, a polyhedral map $K$
is $\{p, q\}$-equivelar if and only if $\widetilde{K}$ is $\{q,
p\}$-equivelar.

Since the degrees of the vertices in an equivelar combinatorial
2-manifold are same, an equivelar combinatorial 2-manifold is also
called degree regular. More explicitly, a connected combinatorial
$2$-manifold $X$ is said to be {\em degree-regular of type $d$}
if each vertex of $X$ has degree $d$. A combinatorial
$2$-manifold is said to be {\em degree-regular} if it is
degree-regular of type $d$ for some $d$. So, trivial examples of
degree-regular combinatorial $2$-manifolds are neighbourly
combinatorial $2$-manifolds.

A {\em combinatorially regular} combinatorial 2-manifold is a
connected combinatorial 2-manifold with a flag-transitive
automorphism group (a {\em flag} is a triple $(u, e, F)$, where
$e$ is an edge of the face $F$ and $u$ is a vertex of $e$). A
connected combinatorial $2$-manifold $X$ is said to be {\em
weakly regular} if the automorphism group of $X$ acts
transitively on $V(X)$. Clearly, a combinatorially regular
combinatorial $2$-manifold is weakly regular and a weakly regular
combinatorial $2$-manifold is degree-regular. If a degree-regular
(respectively, weakly regular) combinatorial 2-manifold $X$
triangulates a space $Y$ then we say $X$ is a {\em
degree-regular} (respectively, {\em weakly regular}) {\em
triangulation} of $Y$.

If the number of $i$-simplices of a $d$-dimensional finite
simplicial complex $X$ is $f_i(X)$ ($0\leq i\leq d$), then the
number $\chi(X) := \sum_{i= 0}^{d}(- 1)^i f_i(X)$ is called the
{\em Euler characteristic} of $X$. Similarly, for a polyhedral
map $K$, the number of vertices $-$ the number of edges $+$ the
number of faces is called the {\em Euler characteristic} of $K$.

For the existence of an $n$-vertex neighbourly combinatorial
2-manifold, $n(n - 1)$ must be divisible by 6, equivalently, $n
\equiv 0$ or 1 mod 3. Ringel and Jungerman (\cite{jr, r}) have
shown that {\em there exist neighbourly combinatorial
$2$-manifolds on $3k$ and $3k + 1$ vertices, for each $k \geq
2$}. In \cite{ab}, Altshuler and Brehm have shown that {\em there
are exactly $19$ neighbourly combinatorial $2$-manifolds on at
most $11$ vertices}. By using computer, Altshuler {\it et al}
(\cite{abs}) have shown that {\em there are exactly $59$
orientable neighbourly combinatorial $2$-manifolds on $12$
vertices}. In \cite{bg}, Bokowski and Guedes de Oliveira have
shown that one of these 59 combinatorial 2-manifolds (namely,
$N^{12}_{54}$) is not geometrically realizable in $\RR^3$. In
\cite{dn}, Datta and Nilakantan have shown that {\em there are
exactly $27$ degree-regular combinatorial $2$-manifolds on at
most $11$ vertices}.

If the degree of each vertex in an $n$-vertex combinatorial
2-manifold $K$ is $d$ then $nd = 2f_1(K) = 3f_2(K)$ and $\chi(K)
= f_0(K) - f_1(K) + f_2(K) = n - \frac{nd}{2} + \frac{nd}{3} =
\frac{n(6 - d)}{6}$. So, if $K$ is weakly regular and $\chi(K)
\geq 0$ then $(n, d) = (4, 3), (6, 4), (6, 5), (12, 5)$ or
$(\chi(K), d) = (0, 6)$. For each $(n, d) \in \{(4, 3), (6, 4),
(6, 5), (12, 5)\}$, there exists a unique combinatorial
2-manifold, namely, the 4-vertex 2-sphere, the boundary of the
octahedron, the 6-vertex real projective plane and the boundary
of the icosahedron (cf. \cite{bd, dn}). These 4 combinatorial
2-manifolds are combinatorially regular. Any degree-regular
triangulation of the torus is weakly regular and for each $n\geq
7$, there exists an $n$-vertex weakly regular triangulation of
the torus. There exists an $n$-vertex degree-regular
triangulation of the Klein bottle if and only if $n$ is a
composite number $\geq 9$ and for $k \geq 2$, there exists a $(4k
+ 2)$-vertex weakly regular triangulation of the Klein bottle
(cf. \cite{dn, du1}). There are infinitely many combinatorially
regular triangulations of the torus whereas there is no such
triangulation of the Klein bottle (cf. \cite{v}). For the
existence of degree-regular of type $d$ combinatorial 2-manifolds
of negative Euler characteristic, $d$ must be at least 7. Since
$\frac{n(6-d)}{6} \neq -1$ for $n > d \geq 7$, there does not
exist any degree-regular combinatorial 2-manifolds of Euler
characteristic $-1$. Schulte and Wills (\cite{sw1, sw2}) have
constructed two combinatorially regular geometrically realized
triangulations of the orientable surface of genus 3 in $\RR^3$.

The orientable surface of genus 2 (the double-torus) has Euler
characteristic $-2$. It is known that for a triangulation of the
double-torus one needs at least 10 vertices. In \cite{l2}, Lutz
has shown that there are exactly 865 distinct 10-vertex
triangulations of the double torus and all of them are
geometrically realizable in $\RR^3$. If $K$ is an $n$-vertex
degree-regular of type $d$ combinatorial 2-manifold of Euler
characteristic $-2$ then $-2 = \frac{n(6 - d)}{6}$ and hence $(n,
d) = (12, 7)$. In \cite{msw}, McMullen {\em et al.} constructed a
weakly regular triangulation (isomorphic to $N_1$ of Example 1)
of the double-torus which has a linear geometric realization in
$\RR^3$. In \cite{l1}, Lutz showed that {\em there are exactly
$3$ orientable weakly regular combinatorial $2$-manifolds of
Euler characteristic $-2$ $($i.e., triangulations of the
double-torus$)$. These $3$ are $N_1$, $N_2$ and $N_3$ of Example
$\ref{g2e1}$}. Lutz also showed that ${\rm Aut}(N_1) \cong {\rm
Aut}(N_2) \cong \ZZ_{12}$ and ${\rm Aut}(N_3) \cong D_{6}$. Here
we prove the following\,:


\begin{theo}$\!\!${\bf .}
Let $N_1, \dots, N_6$ be as in Example $\ref{g2e1}$. Then
\begin{enumerate}
    \item[{\rm (a)}] $N_i \not \cong N_j$ for $1\leq i< j\leq 6$,
    \item[{\rm (b)}] ${\rm Aut}(N_4) \cong \ZZ_{2}$, ${\rm Aut}(N_5)
    \cong D_3$, ${\rm Aut}(N_6) \cong \ZZ_{2}\times \ZZ_{2}$, and
    \item[{\rm (c)}] if $M$ is an orientable degree-regular
    combinatorial $2$-manifold of Euler characteristic $-2$ then
    $M$ is isomorphic  to one of $N_1, \dots, N_6$.
\end{enumerate}
\end{theo}

If there exists an $n$-vertex $\{p, q\}$-equivelar polyhedral map
of Euler characteristic $-2$ then $-2 = n - \frac{nq}{2} +
\frac{nq}{p}$ and $\frac{nq(p - 3)}{2} + \frac{nq}{2} \leq {n
\choose 2}$. So, $(n, p, q) = (12, 3, 7)$ or $(28, 7, 3)$. Since
two polyhedral maps $M$ and $N$ are isomorphic if and only if
their duals $\widetilde{M}$ and $\widetilde{N}$ are isomorphic
(cf. \cite{ms}), from Theorem 1, we get.

\begin{cor}$\!\!${\bf .}
There are exactly $12$ orientable equivelar polyhedral maps of
Euler characteristic $-2$, namely, $N_1, \dots, N_6$,
$\widetilde{N_1}, \dots, \widetilde{N_5}$ and $\widetilde{N_6}$.
\end{cor}

\section{Examples and Preliminaries}

In this section we present those degree-regular combinatorial
2-manifolds and their automorphism groups mentioned in the
previous section. First we give some definitions and notations
which will be used throughout the paper.

A 2-simplex in a 2-dimensional simplicial complex is also said to
be a {\em face}. We denote a face $\{u, v, w\}$ by $uvw$. We also
denote an edge $\{u, v\}$ by $uv$.

A {\em graph} is an one dimensional simplicial complex. The
complete graph on $n$ vertices is denoted by $K_n$.  Disjoint
union of $m$ copies of $K_n$ is denoted by $mK_n$. A graph
without any edge is called a {\em null} graph. An $n$-vertex null
graph is denoted by $\emptyset_{n}$.

If $G$ is a graph and $n \geq 0$ is an integer then we consider
the graph $G_n(G)$ as follows. The vertices of $G_n(G)$ are the
vertices of $G$. Two vertices $u$ and $v$ form an edge in
$G_n(G)$ if the number of common neighbours of $u$ and $v$ is
$n$. Clearly, if $G$ and $H$ are isomorphic then $G_n(G)$ and
$G_n(H)$ are isomorphic for all $n \geq 0$.

A connected finite graph is called a {\em cycle} if the degree of
each vertex is 2. An $n$-cycle is a cycle on $n$ vertices and is
denoted by $C_n$ (or by $C_n(a_1, \dots, a_n)$ if the edges are
$a_1a_2, \dots, a_{n-1}a_n, a_na_1$).  Disjoint union of $m$
copies of $C_n$ is denoted by $mC_n$.

For a simplicial complex $K$, the graph consisting of the edges
and vertices of $K$ is called the {\em edge-graph} of $K$ and is
denoted by ${\rm EG}(K)$. Let $K$ be a simplicial complex with
vertex-set $V(K)$. If $U \subseteq V(K)$ then the {\em induced
subcomplex} of $K$ on $U$ (denoted by $K[U]$) is the subcomplex
whose simplices are those of $K$ which are subsets of $U$.

If $v$ is a vertex of a simplicial complex $X$, then the {\em
link} of $v$ in $X$, denoted by ${\rm lk}_X(v)$ (or ${\rm
lk}(v)$), is the simplicial complex $\{\tau \in X ~ \colon ~ v
\not\in \tau, ~ \{v\}\cup\tau\in X\}$. If $v$ is a vertex of a
simplicial complex $X$, then the {\em star} of $v$ in $X$,
denoted by ${\rm st}_X(v)$ (or ${\rm st}(v)$), is the simplicial
complex $\{\{v\}, \tau, \tau \cup \{v\} ~ \colon ~ \tau \in {\rm
lk}_X(v)\}$. Clearly, a finite simplicial complex $K$ is a
combinatorial $2$-manifold if and only if ${\rm lk}_K(v)$ is a
cycle for each vertex $v$ of $K$.

\begin{eg}$\!\!${\bf .} \label{g2e1}
{\rm Six degree-regular combinatorial $2$-manifolds of Euler
characteristic $-2$.  Each of them is orientable and hence
triangulates the connected sum of two copies of the torus (the
double-torus). The vertex-set of each of these is $\{0, 1, \dots,
11\}$.}

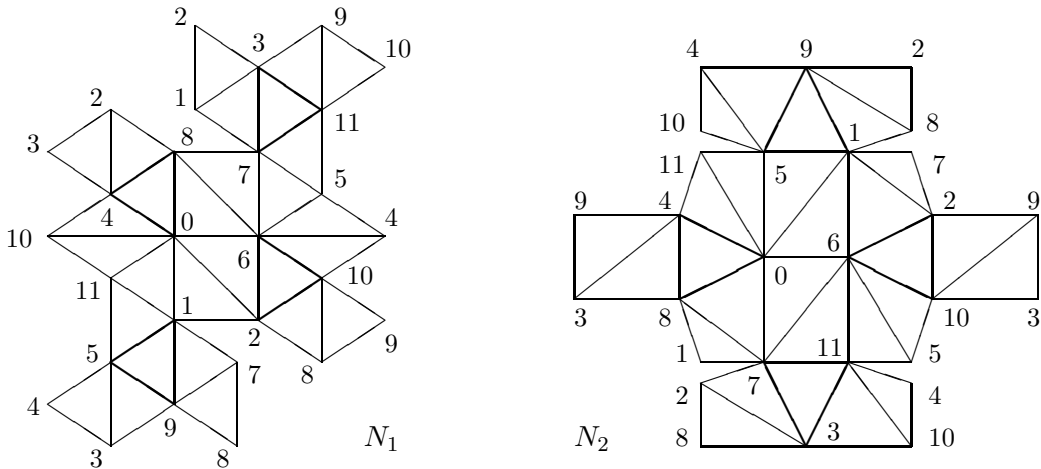
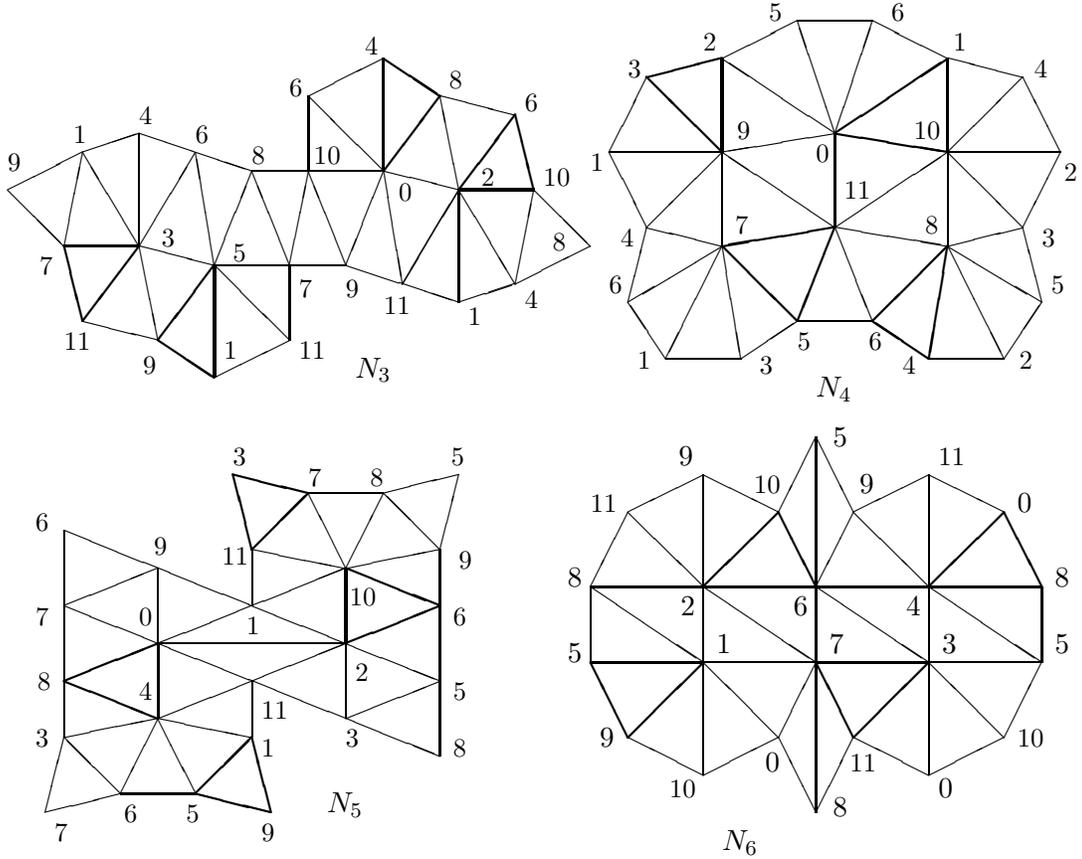
\setlength{\unitlength}{2.8mm}

\begin{picture}(50,24)(0,0)


\thicklines

\put(9,4){\line(0,1){4}} \put(9,12){\line(0,1){4}}
\put(13,8){\line(0,1){4}} \put(13,16){\line(0,1){4}}

\put(6,6){\line(3,2){3}} \put(6,6){\line(3,-2){3}}
\put(6,14){\line(3,2){3}} \put(6,14){\line(3,-2){3}}
\put(13,8){\line(3,2){3}} \put(13,12){\line(3,-2){3}}
\put(13,16){\line(3,2){3}} \put(13,20){\line(3,-2){3}}

\thinlines

\put(9,8){\line(1,0){4}} \put(3,12){\line(1,0){16}}
\put(9,16){\line(1,0){4}}

\put(6,2){\line(0,1){8}} \put(6,14){\line(0,1){4}}
\put(9,8){\line(0,1){4}} \put(12,2){\line(0,1){4}}
\put(10,18){\line(0,1){4}} \put(13,12){\line(0,1){4}}
\put(16,6){\line(0,1){4}} \put(16,14){\line(0,1){8}}

\put(3,4){\line(3,2){3}} \put(3,12){\line(3,2){3}}
\put(3,16){\line(3,2){3}} \put(6,2){\line(3,2){6}}
\put(6,10){\line(3,2){3}} \put(16,6){\line(3,2){3}}
\put(16,10){\line(3,2){3}} \put(13,12){\line(3,2){3}}
\put(10,18){\line(3,2){6}} \put(16,18){\line(3,2){3}}

\put(3,4){\line(3,-2){3}} \put(3,12){\line(3,-2){6}}
\put(3,16){\line(3,-2){3}} \put(6,18){\line(3,-2){3}}
\put(9,4){\line(3,-2){3}} \put(9,8){\line(3,-2){3}}
\put(13,8){\line(3,-2){3}} \put(13,16){\line(3,-2){6}}
\put(16,10){\line(3,-2){3}} \put(16,22){\line(3,-2){3}}
\put(10,22){\line(3,-2){3}} \put(9,12){\line(1,-1){4}}
\put(9,16){\line(1,-1){4}} \put(10,18){\line(3,-2){3}}

\put(2,3.5){\mbox{\small 4}} \put(4.8,6){\mbox{\small 5}}
\put(4.3,9){\mbox{\small 11}} \put(1,11.5){\mbox{\small 10}}
\put(5.5,12.5){\mbox{\small 4}} \put(2,16){\mbox{\small 3}}
\put(5,1){\mbox{\small 3}} \put(8.5,2.5){\mbox{\small 9}}
\put(11,1){\mbox{\small 8}} \put(9.3,8.3){\mbox{\small 1}}
\put(9.3,12.3){\mbox{\small 0}} \put(9.3,16.3){\mbox{\small 8}}
\put(5,18.2){\mbox{\small 2}} \put(12.5,6.8){\mbox{\small 2}}
\put(15,5){\mbox{\small 8}} \put(12.5,5){\mbox{\small 7}}
\put(12,10.5){\mbox{\small 6}} \put(12,14.5){\mbox{\small 7}}
\put(12.7,20.8){\mbox{\small 3}} \put(9,18.2){\mbox{\small 1}}
\put(9,22){\mbox{\small 2}} \put(19,6.5){\mbox{\small 9}}
\put(19,12.5){\mbox{\small 4}} \put(17.2,9.7){\mbox{\small 10}}
\put(16.6,14.3){\mbox{\small 5}} \put(16.6,17){\mbox{\small 11}}
\put(16.6,22){\mbox{\small 9}} \put(19,20.5){\mbox{\small 10}}


\put(18,2){\mbox{$N_{1}$}}



\thicklines

\put(33,9){\line(0,1){4}} \put(45,9){\line(0,1){4}}
\put(37,6){\line(1,0){4}} \put(37,16){\line(1,0){4}}
\put(39,2){\line(-1,2){2}} \put(39,2){\line(1,2){2}}
\put(39,20){\line(-1,-2){2}} \put(39,20){\line(1,-2){2}}
\put(37,11){\line(-2,-1){4}} \put(37,11){\line(-2,1){4}}
\put(41,11){\line(2,-1){4}} \put(41,11){\line(2,1){4}}

\thinlines

\put(28,9){\line(1,0){5}} \put(45,9){\line(1,0){5}}
\put(45,13){\line(1,0){5}} \put(28,13){\line(1,0){5}}
\put(34,2){\line(1,0){10}} \put(34,6){\line(1,0){10}}
\put(34,16){\line(1,0){10}} \put(34,20){\line(1,0){10}}
\put(37,11){\line(1,0){4}}

\put(28,9){\line(0,1){4}} \put(50,9){\line(0,1){4}}
\put(34,2){\line(0,1){3}} \put(44,2){\line(0,1){3}}
\put(34,17){\line(0,1){3}} \put(44,17){\line(0,1){3}}
\put(37,6){\line(0,1){10}} \put(41,6){\line(0,1){10}}
\put(34,5){\line(3,1){3}} \put(34,6){\line(-1,3){1}}
\put(33,13){\line(1,3){1}} \put(37,16){\line(-3,1){3}}
\put(41,16){\line(3,1){3}} \put(44,16){\line(1,-3){1}}
\put(44,6){\line(1,3){1}} \put(41,6){\line(3,-1){3}}
\put(45,9){\line(5,4){5}} \put(28,9){\line(5,4){5}}
\put(39,2){\line(-5,3){5}} \put(44,2){\line(-3,4){3}}
\put(44,17){\line(-5,3){5}} \put(37,16){\line(-3,4){3}}
\put(37,6){\line(-4,3){4}} \put(44,6){\line(-3,5){3}}
\put(45,13){\line(-4,3){4}} \put(37,11){\line(-3,5){3}}
\put(37,6){\line(4,5){4}} \put(37,11){\line(4,5){4}}

\put(33.3,20.5){\mbox{\small 4}} \put(38.7,20.5){\mbox{\small 9}}
\put(44,20.5){\mbox{\small 2}} \put(32,17){\mbox{\small 10}}
\put(41,16.5){\mbox{\small 1}} \put(44.7,17){\mbox{\small 8}}
\put(32,15){\mbox{\small 11}} \put(37.5,14.5){\mbox{\small 5}}
\put(45,15){\mbox{\small 7}}

\put(28,13.3){\mbox{\small 9}} \put(32,13.3){\mbox{\small 4}}
\put(45.5,13.3){\mbox{\small 2}} \put(49.5,13.3){\mbox{\small 9}}
\put(37.5,9.7){\mbox{\small 0}} \put(40,11.3){\mbox{\small 6}}
\put(28,7.7){\mbox{\small 3}} \put(32,7.7){\mbox{\small 8}}
\put(45.5,7.7){\mbox{\small 10}} \put(49.5,7.7){\mbox{\small 3}}
\put(32.8,6){\mbox{\small 1}} \put(44.8,6){\mbox{\small 5}}
\put(39.5,6.3){\mbox{\small 11}} \put(36.2,4.5){\mbox{\small 7}}
\put(32.8,4){\mbox{\small 2}} \put(44.8,4){\mbox{\small 4}}
\put(32.8,2){\mbox{\small 8}} \put(44.8,2){\mbox{\small 10}}
\put(40,2.3){\mbox{\small 3}}

\put(28,2){\mbox{$N_2$}}


\end{picture}

\setlength{\unitlength}{2.5mm}

\begin{picture}(55,22)(0,0)


\thicklines

\put(4,8){\line(1,-4){1}} \put(4,8){\line(1,0){4}}
\put(8,8){\line(-3,-4){3}} \put(12,7){\line(-3,-4){3}}
\put(12,1){\line(0,1){6}} \put(9,3){\line(3,-2){3}}
\put(21,12){\line(0,1){6}} \put(21,12){\line(3,4){3}}
\put(21,18){\line(3,-2){3}} \put(25,11){\line(1,0){4}}
\put(25,11){\line(3,4){3}} \put(28,15){\line(1,-4){1}}

\thinlines

\put(1,11){\line(2,1){4}} \put(5,13){\line(3,1){3}}
\put(8,14){\line(3,-1){6}} \put(14,12){\line(1,0){7}}
\put(17,12){\line(0,1){4}} \put(17,16){\line(2,1){4}}
\put(24,16){\line(4,-1){4}} \put(29,11){\line(1,-1){3}}
\put(1,11){\line(1,-1){3}} \put(5,4){\line(4,-1){4}}
\put(12,1){\line(2,1){4}} \put(16,3){\line(0,1){4}}
\put(12,7){\line(1,0){7}} \put(19,7){\line(3,-1){6}}
\put(25,5){\line(3,1){3}} \put(28,6){\line(2,1){4}}
\put(8,8){\line(3,5){3}} \put(22,6){\line(3,5){3}}
\put(4,8){\line(1,5){1}} \put(8,8){\line(-3,5){3}}
\put(8,8){\line(0,1){6}} \put(8,8){\line(4,-1){4}}
\put(8,8){\line(1,-5){1}} \put(12,7){\line(-1,6){1}}
\put(12,7){\line(2,5){2}} \put(12,7){\line(1,-1){4}}
\put(16,7){\line(-2,5){2}} \put(16,7){\line(1,5){1}}
\put(19,7){\line(-2,5){2}} \put(28,6){\line(1,5){1}}
\put(21,12){\line(-2,-5){2}} \put(21,12){\line(1,-6){1}}
\put(21,12){\line(-1,1){4}} \put(21,12){\line(4,-1){4}}
\put(25,11){\line(-1,5){1}} \put(25,11){\line(3,-5){3}}
\put(25,11){\line(0,-1){6}} \put(25,11){\line(-3,-5){3}}

\put(9.2,8){\mbox{\small 3}} \put(14,12.4){\mbox{\small 8}}
\put(17.3,12.4){\mbox{\small 10}} \put(8,14.5){\mbox{\small 4}}
\put(13,7.2){\mbox{\small 5}} \put(16,16){\mbox{\small 6}}
\put(1,12){\mbox{\small 9}} \put(4.5,13.5){\mbox{\small 1}}
\put(11,13.5){\mbox{\small 6}} \put(20,18.2){\mbox{\small 4}}
\put(24.5,16.4){\mbox{\small 8}} \put(28.5,15.2){\mbox{\small 6}}
\put(29.5,11.2){\mbox{\small 10}} \put(30,7.8){\mbox{\small 8}}
\put(28.5,4.7){\mbox{\small 4}} \put(25.5,3.8){\mbox{\small 1}}
\put(21,4.4){\mbox{\small 11}} \put(19,5.3){\mbox{\small 9}}
\put(16.5,5.3){\mbox{\small 7}} \put(16.5,2){\mbox{\small 11}}
\put(12.5,2){\mbox{\small 1}} \put(8.2,1.2){\mbox{\small 9}}
\put(4,2.5){\mbox{\small 11}} \put(2.7,6.5){\mbox{\small 7}}
\put(21.8,10.4){\mbox{\small 0}} \put(26.2,11.2){\mbox{\small 2}}

\put(19.5,1){\mbox{$N_{3}$}}


\thicklines

\put(39,13){\line(-1,1){4}} \put(39,13){\line(0,1){5}}
\put(35,17){\line(4,1){4}} \put(45,14){\line(6,-1){6}}
\put(45,14){\line(3,2){6}} \put(51,13){\line(0,1){5}}
\put(45,9){\line(-2,-5){2}} \put(45,9){\line(-6,-1){6}}
\put(39,8){\line(1,-1){4}} \put(51,8){\line(-1,-1){4}}
\put(51,8){\line(-1,-6){1}} \put(50,2){\line(-3,2){3}}

\thinlines

\put(36,2){\line(1,0){4}} \put(50,2){\line(1,0){4}}
\put(33,13){\line(1,0){6}} \put(51,13){\line(1,0){6}}
\put(43,4){\line(1,0){4}} \put(43,20){\line(1,0){4}}
\put(51,8){\line(0,1){5}} \put(39,8){\line(0,1){5}}

\put(45,9){\line(0,1){5}} \put(45,9){\line(2,-5){2}}
\put(45,9){\line(6,-1){6}} \put(45,9){\line(3,2){6}}
\put(45,9){\line(-3,2){6}}

\put(45,14){\line(-1,3){2}} \put(45,14){\line(1,3){2}}
\put(45,14){\line(-3,2){6}} \put(45,14){\line(-6,-1){6}}

\put(39,13){\line(-1,-1){4}} \put(51,13){\line(1,1){4}}
\put(51,13){\line(1,-1){4}}

\put(39,8){\line(1,-6){1}} \put(39,8){\line(-1,-2){3}}
\put(39,8){\line(-5,-3){5}} \put(39,8){\line(-4,1){4}}

\put(51,8){\line(1,-2){3}} \put(51,8){\line(5,-3){5}}
\put(51,8){\line(4,1){4}}

\put(36,2){\line(-2,3){2}} \put(34,5){\line(1,4){1}}
\put(33,13){\line(1,-2){2}} \put(33,13){\line(1,2){2}}
\put(39,18){\line(2,1){4}} \put(55,17){\line(-4,1){4}}
\put(51,18){\line(-2,1){4}} \put(57,13){\line(-1,-2){2}}
\put(57,13){\line(-1,2){2}} \put(54,2){\line(2,3){2}}
\put(56,5){\line(-1,4){1}} \put(40,2){\line(3,2){3}}

\put(39.7,8.7){\mbox{\small 7}} \put(45.5,10.5){\mbox{\small 11}}
\put(49.8,8.7){\mbox{\small 8}} \put(39.8,13.7){\mbox{\small 9}}
\put(44,12.5){\mbox{\small 0}} \put(49.2,13.7){\mbox{\small 10}}
\put(32,12){\mbox{\small 1}} \put(33.6,8){\mbox{\small 4}}
\put(33,5.5){\mbox{\small 6}} \put(34.5,1.5){\mbox{\small 1}}
\put(41,1.2){\mbox{\small 3}} \put(43,2.5){\mbox{\small 5}}
\put(46.8,2.5){\mbox{\small 6}} \put(48.6,1.2){\mbox{\small 4}}
\put(54.8,1.5){\mbox{\small 2}} \put(57.2,11.5){\mbox{\small 2}}
\put(56,8){\mbox{\small 3}} \put(56.5,5.5){\mbox{\small 5}}
\put(34,17){\mbox{\small 3}}  \put(38,18.4){\mbox{\small 2}}
\put(41.5,20){\mbox{\small 5}} \put(55.6,17){\mbox{\small 4}}
\put(51.3,18.4){\mbox{\small 1}} \put(48,20){\mbox{\small 6}}

\put(44,0){\mbox{$N_{4}$}}


\end{picture}

\setlength{\unitlength}{2.5mm}

\begin{picture}(55,25)(0,-1)


\thicklines

\put(11,3){\line(1,1){3}} \put(11,3){\line(4,-1){4}}
\put(15,2){\line(-1,4){1}} \put(9,7){\line(0,1){4}}
\put(4,9){\line(5,2){5}} \put(4,9){\line(5,-2){5}}
\put(19,11){\line(0,1){4}}  \put(19,15){\line(5,-2){5}}
\put(19,11){\line(5,2){5}} \put(14,16){\line(-1,4){1}}
\put(14,16){\line(1,1){3}} \put(13,20){\line(4,-1){4}}

\thinlines

\put(7,3){\line(1,0){4}} \put(9,11){\line(1,0){10}}
\put(17,19){\line(1,0){4}} \put(4,6){\line(0,1){11}}
\put(19,7){\line(0,1){4}} \put(24,5){\line(0,1){11}}
\put(14,6){\line(0,1){3}} \put(14,13){\line(0,1){3}}
\put(9,11){\line(0,1){4}}

\put(3,2){\line(4,1){4}} \put(3,2){\line(1,4){1}}
\put(7,3){\line(-1,1){3}} \put(4,13){\line(5,-2){20}}
\put(4,17){\line(5,-2){20}} \put(9,7){\line(5,2){10}}
\put(9,11){\line(5,2){10}}

\put(9,7){\line(5,-1){5}} \put(14,16){\line(5,-1){5}}
\put(11,3){\line(4,-1){4}} \put(4,6){\line(5,1){5}}
\put(19,15){\line(5,1){5}} \put(21,19){\line(4,1){4}}
\put(4,13){\line(5,2){5}} \put(19,7){\line(5,2){5}}
\put(24,16){\line(1,4){1}} \put(24,16){\line(-1,1){3}}
\put(9,7){\line(-1,-2){2}} \put(9,7){\line(1,-2){2}}
\put(19,15){\line(-1,2){2}} \put(19,15){\line(1,2){2}}

\put(3.5,0.5){\mbox{\small 7}} \put(2.5,5.5){\mbox{\small 3}}
\put(2.6,8.5){\mbox{\small 8}} \put(2.5,12){\mbox{\small 7}}
\put(2.5,17){\mbox{\small 6}} \put(8.8,15.7){\mbox{\small 9}}
\put(7.2,1.5){\mbox{\small 6}} \put(10.5,1.5){\mbox{\small 5}}
\put(14.5,0.5){\mbox{\small 9}} \put(14.5,5){\mbox{\small 1}}
\put(14.5,7){\mbox{\small 11}} \put(19.5,9){\mbox{\small 2}}
\put(13.7,11.5){\mbox{\small 1}} \put(12.4,15.2){\mbox{\small
11}} \put(13,20.5){\mbox{\small 3}} \put(8,12){\mbox{\small 0}}
\put(8,8){\mbox{\small 4}} \put(19.2,13){\mbox{\small 10}}
\put(20.3,19.4){\mbox{\small 8}} \put(17,19.4){\mbox{\small 7}}
\put(24.6,20.5){\mbox{\small 5}} \put(24.7,8){\mbox{\small 5}}
\put(24.7,12){\mbox{\small 6}}  \put(25,15){\mbox{\small 9}}
\put(24.7,5){\mbox{\small 8}} \put(19,5.5){\mbox{\small 3}}

\put(18,2){\mbox{$N_{5}$}}



\thicklines

\put(32,10){\line(1,0){6}} \put(44,10){\line(1,0){6}}
\put(38,14){\line(1,0){6}} \put(50,14){\line(1,0){6}}
\put(34,6){\line(1,1){4}} \put(34,6){\line(-1,2){2}}
\put(42,18){\line(-1,-1){4}} \put(42,18){\line(1,-2){2}}
\put(54,18){\line(-1,-1){4}} \put(54,18){\line(1,-2){2}}
\put(46,6){\line(1,1){4}} \put(46,6){\line(-1,2){2}}

\thinlines

\put(38,10){\line(1,0){18}} \put(32,14){\line(1,0){18}}
\put(44,2){\line(0,1){20}} \put(38,4){\line(0,1){16}}
\put(50,4){\line(0,1){16}} \put(32,10){\line(0,1){4}}
\put(56,10){\line(0,1){4}}

\put(38,4){\line(-2,1){4}} \put(38,4){\line(2,1){4}}
\put(42,6){\line(-1,1){4}} \put(42,6){\line(1,2){2}}
\put(44,2){\line(1,2){2}} \put(44,2){\line(-1,2){2}}
\put(50,4){\line(-2,1){4}} \put(50,4){\line(2,1){4}}
\put(54,6){\line(-1,1){4}} \put(54,6){\line(1,2){2}}
\put(34,18){\line(1,-1){4}} \put(34,18){\line(-1,-2){2}}
\put(38,20){\line(-2,-1){4}} \put(38,20){\line(2,-1){4}}
\put(44,22){\line(1,-2){2}} \put(44,22){\line(-1,-2){2}}
\put(46,18){\line(1,-1){4}} \put(46,18){\line(-1,-2){2}}
\put(50,20){\line(-2,-1){4}} \put(50,20){\line(2,-1){4}}
\put(32,14){\line(3,-2){6}} \put(38,14){\line(3,-2){6}}
\put(44,14){\line(3,-2){6}} \put(50,14){\line(3,-2){6}}

\put(30.8,10){\mbox{$5$}} \put(38.7,10.5){\mbox{$1$}}
\put(44.7,10.5){\mbox{$7$}} \put(50.7,10.5){\mbox{$3$}}
\put(56.7,10.5){\mbox{$5$}} \put(30.8,14){\mbox{$8$}}
\put(36.8,12.5){\mbox{$2$}} \put(42.8,12.5){\mbox{$6$}}
\put(48.8,12.5){\mbox{$4$}} \put(56.7,14){\mbox{$8$}}
\put(32.5,5.5){\mbox{$9$}} \put(36.2,2.8){\mbox{\small 10}}
\put(41.3,4.2){\mbox{$0$}} \put(44.9,1.8){\mbox{$8$}}
\put(45.8,4.2){\mbox{\small 11}} \put(50.5,2.8){\mbox{$0$}}
\put(54.7,5.5){\mbox{\small 10}} \put(32,18){\mbox{\small 11}}
\put(36.7,20.5){\mbox{$9$}} \put(40.7,19){\mbox{\small 10}}
\put(44.9,21.5){\mbox{$5$}} \put(46.3,19){\mbox{$9$}}
\put(50.5,20.5){\mbox{\small 11}} \put(54.7,18){\mbox{$0$}}

\put(39,0){\mbox{$N_6$}}

\end{picture}

\end{eg}

\begin{lemma}$\!\!${\bf .} \label{g2l1} If $N_1, \dots, N_6$ are as in
Example $\ref{g2e1}$ then $N_i \not \cong N_j$ for $1\leq i <
j\leq 6$.
\end{lemma}

\noindent {\bf Proof.} If $N_{i} \cong N_{j}$ then ${\rm
EG}(N_{i}) \cong {\rm EG}(N_{j})$ and hence $G_n({\rm EG}(N_{i}))
\cong G_n({\rm EG}(N_{j}))$ (as graphs) for all $n\geq 0$.

Observe that $G_2({\rm EG}(N_{1}))= C_{12}(0, 1, \dots, 11)$,
$G_2({\rm EG}(N_{2})) = \emptyset_{12}$, $G_{2}({\rm EG}(N_{3}))
= C_{12}(0, 11, 2, 1, 4, 3, 6, 5, 8, 7, 10, 9)$,  $G_{2}({\rm
EG}(N_{4})) = (\{0, 1, \dots, 11\}, \{\{9, 11\}, \{10, 11\}\})$,
\newline  $G_{2}({\rm EG}(N_5)) = (\{0, 1, \dots, 11\}, \{\{0, 2\}$,
$\{4, 6\}, \{8, 10\}\})$ and $G_{2}({\rm EG}(N_{6})) = (\{0, 1,
\dots, 11\}$, $\{\{5, 8\}, \{6, 7\}, \{9, 11\}, \{10, 0\}\})$.
So, we need to show only that $N_1 \not\cong N_{3}$.

Observe that $G_{5}({\rm EG}(N_{1})) = G_{5}({\rm EG}(N_{3})) =
C_6(1, 3, 5, 7, 9, 11) \cup C_6(0, 2, 4, 6, 8, 10)$. If $\varphi
\colon N_{1} \to N_{3}$ is an isomorphism then we may assume that
$\varphi(0) = 0$ (since $\ZZ_{12}$ acts vertex-transitively on
$N_{1}$). Since $\varphi$ induces an isomorphism from $G_{5}({\rm
EG}(N_{1}))$ to $G_{5}({\rm EG}(N_{3}))$, $\varphi(\{0, 2, \dots,
10\}) = \{0, 2, \dots, 10\}$ and hence $\varphi(\{1, 3, \dots,
11\}) = \{1, 3, \dots, 11\}$. Then, from ${\rm lk}_{N_{1}}(0)$
and ${\rm lk}_{N_{3}}(0)$, $\varphi(\{1, 11\}) = \{9, 11\}$. If
$\varphi(1) = 11$ then, from ${\rm lk}_{N_{1}}(0)$, ${\rm
lk}_{N_{1}}(1)$, ${\rm lk}_{N_{3}}(0)$ and ${\rm
lk}_{N_{3}}(11)$, $\varphi = (11, 9, 7, 5, 3, 1)(4, 6, 8)$. Now,
$\{2, 3, 4\}$ is a face of $N_1$ but $\varphi(\{2, 3, 4\}) = \{1,
2, 6\}$ is not a face of $N_{3}$, a contradiction. If $\varphi(1)
= 9$ then from ${\rm lk}_{N_{1}}(0)$, ${\rm lk}_{N_{1}}(1)$,
${\rm lk}_{N_{3}}(0)$ and ${\rm lk}_{N_{3}}(9)$, we get $\varphi
= (1, 9, 5, 3, 7)(2, 10)(4, 8)$. Now, $\{5, 7, 11\}$ is a face of
$N_1$ but $\varphi(\{5, 7, 11\}) = \{1, 3, 11\}$ is not a face of
$N_{3}$, a contradiction. This completes the proof. \hfill $\Box$

\bigskip

In \cite{l1}, Lutz showed that ${\rm Aut}(N_1) = {\rm Aut}(N_2) =
\langle \sigma \rangle \cong \ZZ_{12}$ and ${\rm Aut}(N_3) =
\langle \sigma^2, \tau \rangle \cong D_{6}$, where $\sigma = (0,
1, \dots, 11)$ and $\tau = (1, 4)(2, 3)(8, 9)(6, 11)(0, 5)(7,
10)$. Here we prove.

\begin{lemma}$\!\!${\bf .} \label{g2l2}
${\rm Aut}(N_4) = \langle \gamma \rangle \cong \ZZ_{2}$, ${\rm
Aut}(N_5) = \langle \alpha, \beta \rangle \cong D_3$ and ${\rm
Aut}(N_6) = \langle \alpha_1, \alpha_2 \rangle \cong \ZZ_{2}
\times \ZZ_{2}$, where $\alpha = (0, 4, 8)(1, 5, 9)(2, 6, 10)(3,
7, 11)$, $\beta = (0, 2)(4, 10)(1, 11)(3, 9)(5, 7)(6, 8)$,
$\alpha_1 = (1, 2)(3, 4)(5, 6)(7, 8)(9, 10)(11, 0)$, $\alpha_2 =
(1, 3)(2, 4)(5, 7)(6, 8)(9, 11)(10, 0)$ and $\gamma = (1, 2)(3,
4)(5, 6)(7, 8)(9, 10)$.
\end{lemma}

\noindent {\bf Proof.} For $i \geq 0$ and $14 \leq j \leq 6$, let
${\cal E}_{i, j}$ denote the edge-set of $G_i({\rm EG}(N_j))$.
Then ${\cal E}_{2,4} = \{\{9, 11\}, \{11, 10\}$, ${\cal E}_{6, 4}
= \{\{3, 4\}, \{9, 10\}\}$ and ${\cal E}_{5, 4} = \{\{11, 1\}$,
$\{11, 2\}$, $\{11, 3\}$, $\{11, 4\}$, $\{0, 3\}$, $\{0, 4\}$,
$\{0, 7\}$, $\{0, 8\}$, $\{1, 2\}$, $\{7, 8\}$, $\{5, 9\}$, $\{5,
10\}$, $\{6, 9\}$, $\{6, 10\}\}$.

Clearly, $\gamma$ is an automorphism of $N_4$ of order 2. Let $f$
be an automorphism. Then $f$ induces an automorphism (also
denoted by $f$) on $G_i({\rm EG}(N_4))$ for all $i\geq 0$.
Therefore, considering the action of $f$ on ${\cal E}_{2, 4}$,
$f(11) = 11$ and hence, considering the action of $f$ on ${\cal
E}_{5, 4}$, $f(0) = 0$. Then, considering the faces through the
edge $\{0, 11\}$, $f(9) = 9$ or 10. In the first case, $f(10) =
10$. Then, from the links of 11 and 0 in $N_4$, $f = {\rm
Id}_{N_4}$. In the second case, $(\gamma \circ f)(9) = 9$ and
hence, by the previous case, $\gamma \circ f = {\rm Id}_{N_4}$.
Thus, $f = \gamma^{-1} = \gamma$. This implies that ${\rm
Aut}(N_4) = \{{\rm Id}_{N_4}, \gamma\} \cong \ZZ_2$.

\smallskip

Clearly, $\alpha, \beta \in {\rm Aut}(N_5)$. Also, $\alpha^3 =
\beta^2 = {\rm Id}_{N_5}$, $\alpha \circ \beta = \beta \circ
\alpha^2$ and $\alpha^2 \circ \beta = \beta \circ \alpha$. So,
$\langle \alpha, \beta \rangle \cong D_3$.

\smallskip

\noindent {\bf Claim 1.} The identity is the only automorphism of
$N_5$ which fixes the vertex 0.

Observe that $G_5({\rm EG}(N_5)) = C_4(0, 3, 10, 5) \cup C_4(4,
7, 2, 9) \cup C_4(8, 11, 6, 1)$,  ${\cal E}_{2, 5} = \{\{0, 2\}$,
$\{4, 6\}, \{8, 10\}\}$ and ${\cal E}_{6, 5} = \{\{0, 10\}, \{2,
4\}, \{6, 8\}\}$.

If $g$ is an automorphism of $N_5$ then $g$ is also an
automorphism of $G_i({\rm EG}(N_5))$ for all $i \geq 0$. Now, if
$g(0) = 0$ then, from the action of $g$ on ${\cal E}_{2, 5}$,
$g(2) = 2$. These imply, considering the action of $g$ on ${\cal
E}_{6, 5}$ and ${\cal E}_{2, 5}$, $g(10) = 10$, $g(4) = 4$, $g(6)
= 6$ and $g(8) = 8$. Then, from the links of 0 and 2 in $N_5$, $g
= {\rm Id}_{N_5}$. This proves the claim.

Now, let $f \in {\rm Aut}(N_5)$. Then, from the action of $f$ on
${\cal E}_{2, 5}$, $f(0) = 0, 2, 4, 6, 8$ or 10.

If $f(0) = 0$ then, by the claim, $f = {\rm Id}_{N_5}$. If $f(0) =
2$ then $(\beta\circ f)(0) = 0$ and hence, by Claim 1, $\beta\circ
f = {\rm Id}_{N_5}$. This implies that $f= \beta^{-1} = \beta$. If
$f(0) = 6$ then $(\beta \circ \alpha^2 \circ f)(0) = 0$ and
hence, by Claim 1, $\beta \circ \alpha^2 \circ f = {\rm
Id}_{N_5}$. This implies that $f = (\beta \circ \alpha^2)^{-1} =
\alpha \circ \beta$. Similarly, $f(0) = 4$ implies $f = \alpha$,
$f(0) = 8$ implies $f = \alpha^2$ and $f(0) = 10$ implies $f=
\alpha^2 \circ \beta$. Thus, $f \in \langle \alpha, \beta
\rangle$. This implies that ${\rm Aut}(N_5) = \langle \alpha,
\beta \rangle \cong D_3$.

Clearly, $\alpha_1$ and $\alpha_2$ are order 2 automorphisms of
$N_6$. Also, $\alpha_1^2 = \alpha_2^2 = {\rm Id}_{N_6}$ and
$\alpha_1 \circ \alpha_2 = \alpha_2 \circ \alpha_1$. Thus,
$\langle \alpha_1, \alpha_2 \rangle \cong \ZZ_2\times \ZZ_2$.

\smallskip

\noindent {\bf Claim 2.} The identity is the only automorphism of
$N_6$ which fixes the vertex 1.

Observe that ${\cal E}_{2, 6} = \{\{5, 8\}$, $\{6, 7\}$, $\{9,
11\}$, $\{0, 10\}\}$ and ${\cal E}_{3, 6} = \{\{1, 0\}$, $\{1,
9\}$, $\{1, 5\}$, $\{1, 7\}$, $\{3, 5\}$, $\{3, 7\}$, $\{3,
10\}$, $\{3, 11\}$, $\{2, 10\}$, $\{2, 11\}$, $\{2, 6\}$, $\{2,
8\}$, $\{4, 6\}$, $\{4, 8\}$, $\{4, 0\}$, $\{4, 9\}\}$.

If $g$ is an automorphism of $N_6$ then $g$ is also an
automorphism of $G_i({\rm EG}(N_6))$ for all $i \geq 0$.
Therefore, from the action of $g$ on ${\cal E}_{3, 6}$, $g(\{1,
2, 3, 4\}) = \{1, 2, 3, 4\}$. Now, let $g(1) = 1$. Since $\{1,
2\}$ is an edge of $N_6$ and $\{1, 3\}$, $\{1, 4\}$ are non-edges
of $N_6$, $g(2) =2$ and hence $g(\{3, 4\}) = \{3, 4\}$. Since
$127$, $128$, $345$, $346$ are faces of $N_6$, $g(\{7, 8\}) = \{7,
8\}$ and $g(\{5, 6\}) = \{5, 6\}$. Then, from the link of 1,
$g(5) = 5$ and hence $g(8) = 8$, $g(7) = 7$, $g(9) = 9$, $g(10) =
10$, $g(0) = 0$ and $g(6) = 6$. Now, from the link of 5, $g(3) =
3$ and $g(4) = 4$. Therefore, $g = {\rm Id}_{N_6}$. This proves
the claim.

Now, let $f \in {\rm Aut}(N_6)$. Then, from the action of $f$ on
${\cal E}_{3, 6}$, $f(1) = 1, 2, 3$ or 4.

If $f(1) = 1$ then, by Claim 2, $f = {\rm Id}_{N_6}$. If $f(1) =
2$ then $(\alpha_1 \circ f)(1) = 1$ and hence, by Claim 2,
$\alpha_1 \circ f = {\rm Id}_{N_d}$. This implies that $f =
\alpha_1^{-1} = \alpha_1$. Similarly, $f(1) = 3$ implies $f=
\alpha_2$ and $f(1) = 4$ implies $f= \alpha_1\circ \alpha_2$.
Thus, ${\rm Aut}(N_6) = \langle \alpha_1, \alpha_2 \rangle \cong
\ZZ_2\times \ZZ_2$. \hfill $\Box$

\bigskip

For an orientable closed surface $S$, ${\rm Hom}^{+}(S)$ denotes
the group of orientation-preserving homeomorphisms of $S$. We say
that a finite group $G$ {\em acts} ({\em effectively and
orientably}) on $S$ if there is a monomorphism $\varepsilon
\colon G \to {\rm Hom}^{+}(S)$. For an orientable combinatorial
2-manifold $M$, if $\sigma \in {\rm Aut}(M)$ then $\sigma$
induces an homeomorphism $|\sigma| \colon |M| \to |M|$ given by
$|\sigma|(f) := f \circ \sigma$. Let $\iota \, \colon H \to {\rm
Aut}(M)$ be a monomorphism. If $|\iota(\alpha)| \in {\rm
Hom}^{+}(|M|)$ for each $\alpha \in H$ then, we have an action
$\varepsilon \colon H \to {\rm Hom}^{+}(|M|)$, namely,
$\varepsilon(\alpha) := |\iota(\alpha)|$. In this case, we say
that $H$ acts {\em simplicially} on $|M|$. In \cite{b}, Broughton
classified all the finite groups which act on the double-torus.
Among others, $\ZZ_2 \times \ZZ_2$, $\ZZ_6$, $D_3$ and $D_6$ act
on the double-torus.  Using Lemma \ref{g2l2}, we get actions of
$\ZZ_2 \times \ZZ_2$, $\ZZ_6$, $D_3$ and $D_6$ on the
double-torus. More explicitly, we have the following.

\begin{lemma}$\!\!${\bf .} \label{g2l4}
$\ZZ_6$ acts simplicially on $|N_1|$ and on $|N_2|$, $D_6$ acts
simplicially on $|N_3|$, $D_3$ acts simplicially on $|N_5|$, and
$\ZZ_2 \times \ZZ_2$ acts simplicially on $|N_6|$.
\end{lemma}

\noindent {\bf Proof.} Since $N_i$ is orientable, $N_i$ has a
coherent orientation for $1\leq i\leq 6$. Fix a coherent
orientation on $N_i$ (say, the faces in the counterclockwise
direction in the pictures as the positively oriented faces,
namely, $+ \{0, 1, 2\} = \langle 0 1 2 \rangle$ in $N_1$) (see
\cite{c} for definitions and notations).

If $\sigma = (0, 1, \dots, 11)$, then $\sigma \in {\rm Aut}(N_i)$
and $|\sigma| : |N_i| \to |N_i|$ is orientation reversing for
$i=1, 2$ . But, $|\sigma^2| : |N_i| \to |N_i|$ is orientation
preserving. So, $|\sigma^2| \in {\rm Hom}^+(|N_i|)$ for $i = 1,
2$. Since the order of $\sigma^2$ is 6, $\bar{1} \mapsto
\sigma^2$ defines a monomorphism $\iota : \ZZ_6 \to {\rm
Aut}(N_i)$ for $i = 1, 2$. This proves the first statement.

If $\tau = (1, 4)(2, 3)(8, 9)(6, 11)(0, 5)(7, 10)$ then ${\rm
Aut}(N_3) = \langle \sigma^2, \tau \rangle$. Easy to check that
$|\sigma^2|$ and $|\tau|$ are in ${\rm Hom}^+(|N_3|)$. The second
statement follows from the fact that $\langle \sigma^2, \tau
\rangle \cong D_{6}$.

If $\alpha$, $\beta$, $\alpha_1$ and $\alpha_2$ are as in Lemma
\ref{g2l2} then, $|\alpha|, |\beta| \in {\rm Hom}^+(|N_5|)$ and
$|\alpha_1|, |\alpha_2| \in {\rm Hom}^+(|N_6|)$. Since, by Lemma
\ref{g2l2}, ${\rm Aut}(N_5) = \langle \alpha, \beta \rangle \cong
D_3$ and ${\rm Aut}(N_6) = \langle \alpha_1, \alpha_2 \rangle
\cong \ZZ_{2} \times \ZZ_{2}$, the last two statements follow.
\hfill $\Box$


\section{Proof of the theorem}

\begin{lemma}$\!\!${\bf .} \label{g2l3} Let $M$ be a combinatorial
$2$-manifold on $12$ vertices. Let $x$, $y$ and $z$ be three
distinct vertices in $M$. If the degree of each vertex in $M$ is
$7$ then the number of faces in ${\rm st}(x) \cup {\rm st}(y)$ is
$\geq 12$ and the number of faces in ${\rm st}(x) \cup {\rm
st}(y) \cup {\rm st}(z)$ is $\geq 15$.
\end{lemma}

\noindent {\bf Proof\,.} Let $m$ be the number of faces in ${\rm
st}(x) \cup \, {\rm st}(y)$. If $xy$ is an edge then $m = 2 + 5 +
5 = 12$. If $xy$ is not an edge  then $m = 7 + 7 = 14$.

Let $n$ be the number of faces in ${\rm st}(x) \cup {\rm st}(y)
\cup {\rm st}(z)$. If $xyz$ is a face then $n = 1 + 3 + (3 \times
4) = 16$. If $xyz$ is not a face but $xy$, $xz$ and $yz$ are
edges then $n = 3\times 2 + 3 \times 3 = 15$. If one of $xy$,
$xz$ and $yz$ is a non-edge, say $xy$ is a non-edge then the
number of faces in ${\rm st}(x) \cup {\rm st}(y)$ is $14$ and
hence $n \geq 17$. \hfill $\Box$

\bigskip

\noindent {\bf Proof of Theorem 1.} Parts (a) and (b) follow from
Lemmas \ref{g2l1} and \ref{g2l2}.

\smallskip

\noindent Part (c)\,: Let $M$ has $n$ vertices and the degree of
each vertex is $d$. Then the number of edges is $nd/2$ and hence
the number of faces is $nd/3$. Therefore $-2 = \chi(M) = n - nd/2
+ nd/3 =  n - nd/6$. This implies that $nd$ is divisible by 6 and
$6< d < n$ and hence $(n, d) = (12, 7)$.

Let the vertex set of $M$ be $V = \{ 0, 1, \dots, 9, u, v \}$ and
$\psi \colon V \to \{0, 1, \dots, 11\}$ be given by $\psi(i) = i$
for $0 \leq i\leq 9$, $\psi(u) = 10$, $\psi(v) = 11$.

\smallskip

\noindent {\bf Claim 1.} There exist faces $abc$ and $bcd$ such
that $ad$ is a non-edge of $M$ and not all edges $ab$, $ac$,
$bd$, $cd$ are in $G_2({\rm EG}(M))$.

Assume without loss of generality, that ${\rm lk}(0) = C_7(1, 2,
\dots, 7)$. Then, there exists $i \in \{1, \dots, 7 \}$ such that
$\{i, i + 1, j \}$ is not a face in $M$ for all $j \in \{ 1,
\dots, 7 \}$, where the addition is modulo 7. (Otherwise,  $\{ i,
i+1, i+3\}$, $\{i, i+1, i+4\}$ or $\{i, i+1, i+5\}$ is a face for
each $i \in \{ 1, \dots, 7 \}$,  then we get 14 faces through the
vertices $0, 1, \dots, 7$. Hence the number of faces through the
remaining four vertices is 14. This is not possible by Lemma
\ref{g2l3}). So, we can assume without loss of generality that
$128$ is a face.

If one of $01$, $02$, $18$, $28$ is not in $G_2({\rm EG}(M))$ then
$012$ and $128$ are as required. So, assume that $01$, $02$,
$18$, $28$ are in $G_2({\rm EG}(M))$. Since 2 and 7 are in the
links of 0 and 1, we may assume that ${\rm lk}(1) = C_7(7, 0, 2,
8, 9, u, v)$. Similarly, since 1 and 3 are in the links of 0 and
2, the link of 2 must be of the form $C_7(8, 1, 0, 3, p, q, r)$,
where $\{p, q, r\} = \{9, u, v\}$. Then $12$ is in $G_5({\rm
EG}(M))$ and $27$ is a non-edge of $M$. These arguments imply that
$012$ and $017$ are as required. This proves the claim.

\smallskip

By Claim 1, we may assume that $012$, $128$ are faces, where $08$
is a non-edge of $M$ and $01$ is not in $G_2({\rm EG}(M))$. Also,
assume that ${\rm lk}(0) = C_7(1, 2, \dots, 7)$. Then, ${\rm
lk}(1) = C_7(8, 2, 0, 7, u_1, u_2, u_3)$ for some $u_1, u_2, u_3
\in V$ and $\{u_1, u_2, u_3\} \neq \{9, u, v\}$. Since $M$ is
orientable, $(u_1, u_2), (u_2, u_3)$ $\notin \{(3, 4), (4, 5), (5,
6)\}$.

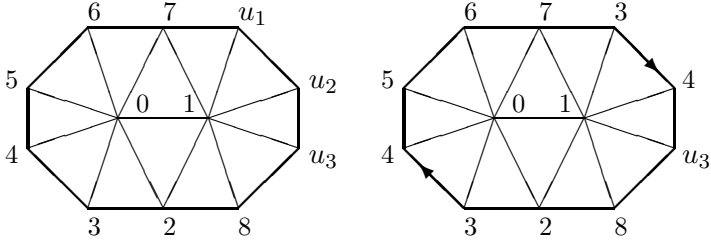
\setlength{\unitlength}{2mm}

\begin{picture}(45,18)(0,-1)


\thicklines

\put(2,6){\line(0,1){4}} \put(2,10){\line(1,1){4}}
\put(6,14){\line(1,0){10}} \put(6,2){\line(1,0){10}}
\put(16,2){\line(1,1){4}} \put(20,6){\line(0,1){4}}
\put(6,2){\line(-1,1){4}} \put(16,14){\line(1,-1){4}}

\thinlines

\put(8,8){\line(1,-2){3}} \put(8,8){\line(-1,-3){2}}
\put(8,8){\line(-3,-1){6}} \put(8,8){\line(-3,1){6}}
\put(8,8){\line(-1,3){2}} \put(8,8){\line(1,2){3}}
\put(8,8){\line(1,0){6}}

\put(14,8){\line(-1,-2){3}} \put(14,8){\line(1,-3){2}}
\put(14,8){\line(3,-1){6}} \put(14,8){\line(3,1){6}}
\put(14,8){\line(1,3){2}} \put(14,8){\line(-1,2){3}}

\put(9.2,8.4){\mbox{\small 0}} \put(12.3,8.4){\mbox{\small 1}}
\put(6,0.2){\mbox{\small 3}} \put(11,0.2){\mbox{\small 2}}
\put(16,0.2){\mbox{\small 8}} \put(6,14.5){\mbox{\small 6}}
\put(11,14.5){\mbox{\small 7}} \put(16,14.5){\mbox{$u_1$}}
\put(0.5,5){\mbox{\small $4$}} \put(0.5,10){\mbox{\small $5$}}
\put(20.7,5){\mbox{$u_3$}} \put(20.7,10){\mbox{$u_2$}}


\thicklines

\put(27,6){\line(0,1){4}} \put(27,10){\line(1,1){4}}
\put(31,14){\line(1,0){10}} \put(31,2){\line(1,0){10}}
\put(41,2){\line(1,1){4}} \put(45,6){\line(0,1){4}}
\put(31,2){\line(-1,1){4}} \put(41,14){\line(1,-1){4}}
\put(31,2){\vector(-1,1){3}} \put(41,14){\vector(1,-1){3}}

\thinlines

\put(33,8){\line(1,-2){3}} \put(33,8){\line(-1,-3){2}}
\put(33,8){\line(-3,-1){6}} \put(33,8){\line(-3,1){6}}
\put(33,8){\line(-1,3){2}} \put(33,8){\line(1,2){3}}
\put(33,8){\line(1,0){6}}

\put(39,8){\line(-1,-2){3}} \put(39,8){\line(1,-3){2}}
\put(39,8){\line(3,-1){6}} \put(39,8){\line(3,1){6}}
\put(39,8){\line(1,3){2}} \put(39,8){\line(-1,2){3}}

\put(34.2,8.4){\mbox{\small 0}} \put(37.3,8.4){\mbox{\small 1}}
\put(31,0.2){\mbox{\small 3}} \put(36,0.2){\mbox{\small 2}}
\put(41,0.2){\mbox{\small 8}} \put(31,14.5){\mbox{\small 6}}
\put(36,14.5){\mbox{\small 7}} \put(41,14.5){\mbox{\small 3}}
\put(25.5,5){\mbox{\small $4$}} \put(25.5,10){\mbox{\small $5$}}
\put(45.5,5){\mbox{$u_3$}} \put(45.5,10){\mbox{\small $4$}}

\end{picture}

Now, it is easy to see that (up to an isomorphism) $(u_1, u_2,
u_3) \in \{(3, 6, 9)$, $(3, 9, 5)$, $(3, 9, 6)$, $(3, 9, u)$,
$(4, 3, 9)$,  $(4, 6, 5)$,  $(4, 9, 6)$, $(4, 9, u)$, $(5, 4,
9)$, $(5, 9, u)$, $(9, 4, u)$, $(9, 5, 4)$, $(9, 5, u)\} \bigcup$
$\{(3, 5, 4)$, $(3, 5, 9)$, $(3, 6, 4)$, $(3, 6, 5)$, $(4, 3,
5)$, $(4, 3, 6)$, $(4, 6, 3)$, $(4, 6, 9)$, $(4, 9, 5)$, $(5, 3,
6)$, $(5, 3, 9)$, $(5, 9, 3)$, $(5, 9, 4)$, $(9, 3, 6)$, $(9, 6,
3)\} \bigcup$ $\{(3, 9, 4)$, $(4, 9, 3)$, $(5, 4, 6)$, $(5, 9,
6)$, $(9, 3, 5)$, $(9, 3, u)$, $(9, 4, 3)$, $(9, 4, 6)$, $(9, 5,
3)$, $(9, 6, 4)$, $(9, 6, 5)$, $(9, 6, u)$, $(9, u, 3)$, $(9, u,
4)$, $(9, u, 5)$, $(9, u, 6)\} = A \cup B \cup C$, say.

\smallskip

\noindent {\bf Claim 2.} A combinatorial 2-manifold corresponding
to any $(u_1, u_2, u_3)\in C$ is isomorphic to a  combinatorial
2-manifold corresponding to some $(u_1, u_2, u_3) \in A \cup B$.

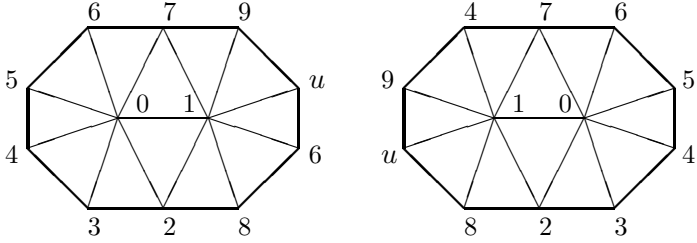
\setlength{\unitlength}{2mm}

\begin{picture}(45,21)(0,-4)


\thicklines

\put(2,6){\line(0,1){4}} \put(2,10){\line(1,1){4}}
\put(6,14){\line(1,0){10}} \put(6,2){\line(1,0){10}}
\put(16,2){\line(1,1){4}} \put(20,6){\line(0,1){4}}
\put(6,2){\line(-1,1){4}} \put(16,14){\line(1,-1){4}}

\thinlines

\put(8,8){\line(1,-2){3}} \put(8,8){\line(-1,-3){2}}
\put(8,8){\line(-3,-1){6}} \put(8,8){\line(-3,1){6}}
\put(8,8){\line(-1,3){2}} \put(8,8){\line(1,2){3}}
\put(8,8){\line(1,0){6}}

\put(14,8){\line(-1,-2){3}} \put(14,8){\line(1,-3){2}}
\put(14,8){\line(3,-1){6}} \put(14,8){\line(3,1){6}}
\put(14,8){\line(1,3){2}} \put(14,8){\line(-1,2){3}}

\put(9.2,8.4){\mbox{\small 0}} \put(12.3,8.4){\mbox{\small 1}}
\put(6,0.2){\mbox{\small 3}} \put(11,0.2){\mbox{\small 2}}
\put(16,0.2){\mbox{\small 8}} \put(6,14.5){\mbox{\small 6}}
\put(11,14.5){\mbox{\small 7}} \put(16,14.5){\mbox{\small $9$}}
\put(0.5,5){\mbox{\small 4}} \put(0.5,10){\mbox{\small $5$}}
\put(20.7,5){\mbox{\small 6}} \put(20.7,10){\mbox{$u$}}
\put(7,-2.5){\mbox{$(0, 1)(3, 8)(4, u, 5, 9, 6)\colon (9, u,
6)\cong (4, 9, u)$}} 


\thicklines

\put(27,6){\line(0,1){4}} \put(27,10){\line(1,1){4}}
\put(31,14){\line(1,0){10}} \put(31,2){\line(1,0){10}}
\put(41,2){\line(1,1){4}} \put(45,6){\line(0,1){4}}
\put(31,2){\line(-1,1){4}} \put(41,14){\line(1,-1){4}}

\thinlines

\put(33,8){\line(1,-2){3}} \put(33,8){\line(-1,-3){2}}
\put(33,8){\line(-3,-1){6}} \put(33,8){\line(-3,1){6}}
\put(33,8){\line(-1,3){2}} \put(33,8){\line(1,2){3}}
\put(33,8){\line(1,0){6}}

\put(39,8){\line(-1,-2){3}} \put(39,8){\line(1,-3){2}}
\put(39,8){\line(3,-1){6}} \put(39,8){\line(3,1){6}}
\put(39,8){\line(1,3){2}} \put(39,8){\line(-1,2){3}}

\put(34.2,8.4){\mbox{\small 1}} \put(37.3,8.4){\mbox{\small 0}}
\put(31,0.2){\mbox{\small 8}} \put(36,0.2){\mbox{\small 2}}
\put(41,0.2){\mbox{\small 3}} \put(31,14.5){\mbox{\small 4}}
\put(36,14.5){\mbox{\small 7}} \put(41,14.5){\mbox{\small 6}}
\put(25.5,5){\mbox{$u$}} \put(25.5,10){\mbox{\small 9}}
\put(45.5,5){\mbox{\small 4}} \put(45.5,10){\mbox{\small 5}}

\end{picture}

Observe that the case $(u_1, u_2, u_3) = (9, u, 6)$ is isomorphic
to the case $(u_1, u_2, u_3) = (4, 9, u)$ by the map $(0, 1)(3,
8)(4, u, 5, 9, 6)$. We denote this by $(0, 1)(3, 8)(4, u, 5, 9,
6) :$ $(9, u, 6)\cong (4, 9, u)$. With this notation we have\,:
\newline
  $(0, 1)(3, 8)(4, u, 5)(6, 9) \, \colon \, (9, u, 5) \cong (9, 4,
  u)$, \hfill $(0, 1)(3, 8)(4, u)(5, 9, 6) \, \colon \,
  (9, 6, u) \cong (5, 9, u)$, \newline
  $(0, 1)(2, 7)(3, 9)(6, 8)(4, 5, u) \, \colon \, (9, u, 4)
  \cong (9, 5, u)$, \hfill $(0, 1)(3, 8)(9, 6, 5) \, \colon \, (9, 6, 4)
  \cong (5, 9, 4)$, \newline
  $(0, 1)(2, 7)(3, 5, u, 4, 9)(6, 8) \, \colon \, (9, u, 3)
  \cong (5, 9, u)$, \hfill $(0, 1)(3, 8)(4, 9, 6, 5) \, \colon \, (9, 6,
  5) \cong (5, 4,  9)$, \newline
  $(0, 1)(2, 7)(3, 5, 4, 9)(6, 8) \, \colon \, (9, 5, 3) \cong (5, 9,
  4)$, \hfill  $(0, 1)(3, 8)(4, 5, 9, 6) \, \colon \, (9, 4, 6)
  \cong (4, 9,  5)$, \newline
  $(0, 1)(2, 7)(6, 8)(3, 4, 9) \, \colon \,(9, 3, 5)
  \cong (4, 9, 5)$, \hfill $(0, 1)(2, 7)(6, 8)(3, 5, 9) \, \colon
  \, (9, 4, 3) \cong (5, 4,  9)$,   \newline
  $(0, 1)(2, 7)(6, 8)(5, u)(3, 4, 9) \, \colon \, (9, 3, u) \cong
  (4, 9, u)$, \hfill $(0, 1)(3, 8)(4, 5, 6) \, \colon \, (5, 4, 6)
  \cong (4, 6, 5)$,  \newline
  $(0, 1)(2, 7)(3, 5, 9, 4)(6, 8) \, \colon \, (4, 9, 3) \cong
  (5, 3,  9)$, \hfill $(0, 1)(3, 8)(4, 9, 5, 6) \, \colon \, (5, 9, 6)
  \cong (4, 6, 9)$,  \newline
  $(0, 1)(2, 7)(6, 8)(4, 5, 9) \, \colon \,(3, 9, 4) \cong (3, 5,
  9)$. \newline
  This proves the claim.

\medskip

\noindent {\bf Claim 3.} There is no orientable combinatorial
2-manifold corresponding to any $(u_1, u_2, u_3)$ in $B$.

\smallskip

If $(u_1, u_2) = (5, 3)$ then, considering ${\rm lk}(5)$ and
orientability of $M$, we see that, $345$ or $567$ is a face. This
is not possible. So, $(u_1, u_2, u_3) \neq (5, 3, 6)$ or $(5, 3,
9)$.

\smallskip

If $(u_1, u_2, u_3) = (3, 5, 4)$ then, considering ${\rm lk}(3)$,
we observe that $347$ or $235$ is a face. Since the set of known
faces remain invariant under the map $(0, 1)(2, 7)(4$, $5)(6,
8)$, we may assume that $347$ is a face. Then, by using Lemma
\ref{g2l3}, we get ${\rm lk}(3) = C_7(5, 1, 7, 4, 0, 2, 9)$, ${\rm
lk}(4) = C_7(8, 1, 5, 0, 3, 7, u)$ and ${\rm lk}(7) = C_7(6, 0,
1, 3, 4, u, v)$. These imply that ${\rm lk}(5) = C_7(9, 3, 1, 4,
0, 6, x)$, for some $x \in V$. It is easy to check that $x \neq
v$. This implies $\deg(v) < 5$. So, $(u_1, u_2, u_3) \neq (3, 5,
4)$.

\smallskip

If $(u_1, u_2, u_3) = (3, 6, 4)$ then $03, 05, 16$ and $18$ are
edges in ${\rm lk}(4)$. Since $346$, $456$ and $458$ can not be
faces, it follows that ${\rm lk}(4) = C_7(5, 0, 3, 8, 1, 6, x)$,
for some $x \in V$. By using Lemma \ref{g2l3}, it is easy to see
that $x = 9$. Then $V({\rm lk}(6)) = \{0, 1, 3, 4, 5, 7, 9\}$ and
$V({\rm lk}(3)) = \{0, 1, 2, 4, 6, 7, 8\}$. These imply $0, 1, 3,
4, 6 \not \in {\rm lk}(u) \cup {\rm lk}(v)$. This is not
possible. Thus, $(u_1, u_2, u_3) \neq (3, 6, 4)$.

\smallskip

If $(u_1, u_2, u_3) = (4, 6, 3)$ then ${\rm lk}(1) = C_7(8, 2, 0,
7, 4, 6, 3)$. Now, considering ${\rm lk}(6)$ we see that $356$ is
a face. This implies, considering ${\rm lk}(3)$, that $235$ a
face. This is not possible, since $M$ is orientable. So, $(u_1,
u_2, u_3) \neq (4, 6, 3)$.

\smallskip

If $(u_1, u_2, u_3) = (5, 9, 3)$ then, considering ${\rm lk}(5)$,
we see that $459 \in M$. Now, considering ${\rm lk}(3)$, we see
that $349$ is a face. This implies $C_4(9, 3, 0, 5) \subseteq
{\rm lk}(4)$. So, $(u_1, u_2, u_3) \neq (5, 9, 3)$.

\smallskip

If $(u_1, u_2, u_3) = (5, 9, 4)$ then, considering ${\rm lk}(5)$
we see that $459$ is a face. Then, $C_3(1, 5, 4)$ $\subseteq {\rm
lk}(9)$. So, $(u_1, u_2, u_3) \neq (5, 9, 4)$.

\smallskip

If $(u_1, u_2, u_3) = (9, 6, 3)$ then, ${\rm lk}(3) = C_7(2, 0,
4, 6, 1, 8, d)$, for some $d \in V$. Now, considering ${\rm
lk}(6)$, we see that $467$ is a face. Then, $M$ is
non-orientable. So, $(u_1, u_2, u_3) \neq (9, 6, 3)$.

By similar arguments, we get $(u_1, u_2, u_3) \neq (3, 5, 9)$,
$(3, 6, 5)$, $(4, 3, 5)$, $(4, 3, 6)$, $ (4, 6, 9)$, $(4, 9, 5)$
and $(9, 3, 6)$.  This proves Claim 3.

\smallskip

In view of Claims 2 and 3, we may assume that $(u_1, u_2, u_3) \in
\{(3, 6, 9)$, $(3, 9, 5)$, $(3, 9, 6)$, $(3, 9, u)$,  $(4, 3,
9)$, $(4, 6, 5)$,  $(4, 9, 6)$, $(4, 9, u)$, $(5, 4, 9)$, $(5, 9,
u)$, $(9, 4, u)$, $(9, 5, 4)$, $(9, 5, u)\}$.

\medskip

\noindent {\bf Case 1.} $(u_1, u_2, u_3) = (3, 6, 9)$, i.e, ${\rm
lk}(1) = C_7(8, 2, 0, 7, 3, 6, 9)$. Now, $05, 07, 13$ and $19$ are
edges in ${\rm lk}(6)$. Since $356$, $367$ or $679$ can not be a
face, it is easy to see that ${\rm lk}(6) = C_7(7, 0, 5, 9, 1, 3,
x)$, where $x = 2, 8$ or $u$.

If $x = 8$ then $V({\rm lk}(3)) = \{0, 1, 2, 4, 6, 7, 8\}$. This
implies that the number of faces through $\{0, 1, \dots, 9\}$ is
$17$. This is not possible by Lemma \ref{g2l3}.

If $x = u$ then $02, 04, 16, 17$ and $6u$ are edges in ${\rm
lk}(3)$. Since $M$ is orientable, ${\rm lk}(3) = C_7(4, 0, 2, u,
6, 1, 7)$. So, $3 \not \in {\rm lk}(v)$ and hence $7 \in {\rm
lk}(v)$. Then ${\rm lk}(7) = C_7(4, 3, 1, 0, 6, u, v)$. This
implies that ${\rm lk}(2) = C_7(u, 3, 0, 1, 8, v, y)$, where, it
is easy to see that $y = 4, 5$ or $9$. If $y = 4$ or 5 then
$\deg(8) < 7$. If $y = 9$ then, considering ${\rm lk}(9)$ and
${\rm lk}(u)$, we get $C_5(1, 2, v, u, 9) \subseteq {\rm lk}(8)$.
So, $x \neq u$. Thus, $x= 2$, i.e., ${\rm lk}(6) = C_7(7, 0, 5,
9, 1, 3, 2)$.

Now, by using Lemma \ref{g2l3}, ${\rm lk}(3) = C_7(4, 0, 2, 6, 1,
7, u)$. Since $0, 1, 3, 6 \not\in {\rm lk}(v)$, $7 \in {\rm
lk}(v)$. Thus, ${\rm lk}(7) = C_7(2, 6, 0, 1, 3, u, v)$. This
implies ${\rm lk}(2) = C_7(8, 1, 0, 3, 6, 7, v)$. Now,
considering the links of $u$, $4$ and $5$, we see that $458$ or
$45v$ is a face. Thus, ${\rm lk}(4) = (u, 3, 0, 5, v, z, w)$ or
$C_7(u, 3, 0, 5, 8, z, w)$ for some $z, w \in V$.

\smallskip

\noindent {\bf Subcase 1.1.} In the first case, considering ${\rm
lk}(8)$, we get $(z, w) = (9, 8)$. So, ${\rm lk}(4) = C_7(u, 3,
0, 5, v$, $9, 8)$. Now, completing successively, we get ${\rm
lk}(8) = C_7(2, 1, 9, 4, u, 5, v)$ and ${\rm lk}(v) = C_7(2, 7$,
$u, 9, 4, 5, 8)$. Then $M \cong N_{1}$ by the map $\psi \circ (1,
8, 7, 4)(2, 6, u, 3)(5, v)$.
\smallskip

\noindent {\bf Subcase 1.2.} In the second case, it is easy to see
that $(z, w) = (v, 9)$. Thus, ${\rm lk}(4) = C_7 (u, 3, 0, 5, 8,
v, 9)$. Now, completing successively, we get ${\rm lk}(8) =
C_7(9, 1, 2, v, 4, 5, u)$, ${\rm lk}(9) = C_7(8, 1, 6, 5, v, 4,
u)$ and ${\rm lk}(5) = C_7(9, 6, 0, 4, 8, u, v)$. Here $M \cong
N_{3}$, by the map $\psi \circ (0, 4, 1, 8, 7)(2, u, v, 9, 5, 3)$.

\medskip

\noindent {\bf Case 2.} $(u_1, u_2, u_3) = (3, 9, 5)$, i.e.,
${\rm lk}(1) = C_7(8, 2, 0, 7, 3, 9, 5)$. Now, $02, 04, 17$ and
$19$ are edges in ${\rm lk}(3)$. So, ${\rm lk}(3) = C_7(4, 0, 2,
9, 1, 7, x)$ or $C_7(2, 0, 4, 7, 1, 9, x)$, for some $x \in V$.
The second case is isomorphic to the first case, by the map $(0,
1)(2, 7)(4, 9)(6, 8)$. So, we may assume ${\rm lk}(3) = C_7(4, 0,
2, 9, 1, 7, x)$ for some $x \in V$. It is easy to see that $x =
8$ or $u$.

If $x = 8$ then, ${\rm lk}(3) = C_7(4, 0, 2, 9, 1, 7, 8)$. Now,
${\rm lk}(2) = C_7(9, 3, 0, 1, 8, y, z)$, for some $y, z \in V$.
If $\{y, z\} \cap \{u, v\} = \emptyset$, then we get $18$ faces
outside ${\rm st}(u) \cup {\rm st}(v)$. This is  not possible by
Lemma \ref{g2l3}. So, up to an isomorphism, $y = u$ or $z = u$.
If $y = u$ then $z = v$, else the vertices $0, 1, 2, 3, 8 \not \in
{\rm lk}(v)$. So, ${\rm lk}(2) = C_7(9, 3, 0, 1, 8, u, v)$. This
implies that ${\rm lk}(8) = C_7(3, 4, u, 2, 1, 5, 7)$ and $V({\rm
lk}(5)) = \{0, 1, 4, 6, 7, 8, 9\}$. Then $0, 1, 3, 5, 8 \not \in
{\rm lk}(v)$. This is not possible. If $z = u$ then $y = 4$ or
$v$. If $y = 4$ then ${\rm lk}(2) = C_7(9, 3, 0, 1, 8, 4, u)$.
Now, it follows that ${\rm lk}(8) = C_7(3, 4, 2, 1, 5, v, 7)$ and
${\rm lk}(4) = C_7(u, 2, 8, 3, 0, 5, v)$. This implies ${\rm
lk}(5) = C_7(9, 1, 8, v, 4, 0, 6)$. Hence ${\rm lk}(9) = C_7(6,
5, 1, 3, 2, u, v)$. But then $C_4(v, 9, 2, 4) \subseteq {\rm
lk}(u)$.  If $y = v$ then, ${\rm lk}(2) = C_7(9, 3, 0, 1, 8, v,
u)$. It is easy to see that ${\rm lk}(8) = C_7(3, 4, v, 2, 1$,
$5, 7)$. So, $V({\rm lk}(5)) = \{0, 1, 4, 6, 7, 8, 9\}$. Then $0,
1, 3, 5, 8 \not \in {\rm lk}(u)$. This is not possible. So, $x
\neq 8$. Thus, $x= u$, i.e., ${\rm lk}(3) = C_7(4, 0, 2, 9, 1, 7,
u)$.

Now ${\rm lk}(2) = C_7(9, 3, 0, 1, 8, b, c)$ for some $b, c \in
V$. It is easy to see that, $(b, c) = (4, 6)$, $(4, u)$, $(4, v)$,
$(6, u)$, $(6, v)$, $(u, 6)$, $(u, v)$, $(v, 4)$, $(v, 6)$ or
$(v, u)$.

If $(b, c) = (4, 6)$ then $0, 1, 2, 3, 4 \not \in {\rm lk}(v)$
and hence $\deg(v) \leq 6$. So, $(b, c) \neq (4, 6)$. If $(b, c)
= (4, v)$ then ${\rm lk}(4) = C_7(2, 8, u, 3, 0, 5, v)$. These
imply that ${\rm lk}(5) = C_7(1, 8, v, 4, 0, 6, 9)$ and $8 \in
{\rm lk}(6)$. Then ${\rm lk}(8) = C_7(6, v, 5, 1, 2, 4, u)$,
${\rm lk}(6) = C_7(8, u, 9, 5, 0, 7, v)$ and ${\rm lk} (9) =
C_7(u, 6, 5, 1, 3, 2, v)$. This implies $\deg(v) = 8$. So, $(b, c)
\neq (4, v)$. If $(b, c) = (v, 4)$ then ${\rm lk}(4) = C_7(2, 9,
u, 3, 0, 5$, $v)$, ${\rm lk}(5) = C_7(1, 9, v, 4, 0, 6, 8)$ and
${\rm lk}(9) = C_7(u, 4, 2, 3, 1, 5, v)$. This implies that $v
\in {\rm lk}(7)$ and hence ${\rm lk}(v) = C_7(u, 9, 5, 4, 2, 8,
7)$. Then ${\rm lk}(7) = C_7(8, v, u, 3$, $1, 0, 6)$ and hence
$C_4(8, 5, 0, 7) \subseteq {\rm lk}(6)$. So, $(b, c) \neq (v,
4)$. By similar arguments, we see that $(b, c) \neq (4, u), (6,
u), (u, 6)$, $(v, 6)$ or $(v, u)$.  Thus, $(b, c) = (6, v)$ or
$(u, v)$.

\smallskip

\noindent {\bf Subcase 2.1.} $(b, c) = (6, v)$, i.e., ${\rm lk}(2)
= C_7(9, 3, 0, 1, 8, 6, v)$. Then, $05, 07, 28$ and $2v$ are edges
in ${\rm lk}(6)$. Since $568$ can not be a face, ${\rm lk}(6) =
C_7(7, 0, 5, v, 2, 8, d)$ or $C_7(v, 2, 8, 7, 0, 5, d)$, for some
$d \in V$. It is easy to see that ${\rm lk}(6) = C_7(7, 0, 5, v,
2, 8$, $4)$ or $C_7(v, 2, 8, 7, 0, 5, u)$.

If ${\rm lk}(5) = C_7(9, 1, 8, 4, 0, 6, v)$ then $C_5(v, 5, 1, 3,
2) \subseteq {\rm lk}(9)$. This is not possible. So, ${\rm lk}(6)
= C_7(v, 2, 8, 7, 0, 5, u)$. Completing successively, we get ${\rm
lk}(7) = C_7(6, 0, 1, 3, u, v, 8)$, ${\rm lk}(u) = C_7(5, 6, v,
7, 3, 4, 9)$, ${\rm lk}(5) = C_7(8, 1, 9, u, 6, 0, 4)$, ${\rm
lk}(8) = C_7(4, 5, 1, 2, 6, 7, v)$ and ${\rm lk}(4) = C_7(u, 9,
v, 8, 5, 0, 3)$. Then $M \cong N_{4}$ by the map $\psi \circ (0,
9, 6, 2, 4, v, 8, u, 5)(3, 7)$.

\smallskip

\noindent {\bf Subcase 2.2.} $(b, c) = (u, v)$, i.e., ${\rm lk}(2)
= C_7(9, 3, 0, 1, 8, u, v)$. Now, it is easy to see that ${\rm
lk}(9) = C_7(5, 1, 3, 2, v$, $y, x)$, where $(x, y) =(4, 8)$,
$(4, 6)$, $(6, 8)$ or $(6, 4)$.

If $(x, y) = (4, 8)$ then, ${\rm lk}(8) = C_7(4, u, 2, 1, 5, v,
9)$. This implies that $C_6(u, 8, 9, 5, 0, 3) \subseteq {\rm
lk}(4)$. If $(x, y) = (4, 6)$ then, it is easy to see that ${\rm
lk}(5) = C_7(8, 1, 9, 4, 0, 6, v)$. This implies that ${\rm
lk}(6) = C_7(4, 9, v, 5, 0, 7, 8)$. Then $\deg(8) \geq 8$. If
$(x, y) = (6, 8)$ then, we get ${\rm lk}(8) = C_7(6, 9, v, u, 2,
1, 5)$. This implies $C_4(6, 8, 1, 9) \subseteq {\rm lk}(5)$. So,
$(x, y) = (6, 4)$.

Now, completing successively, we get ${\rm lk}(4) = C_7(v, 5, 0,
3$, $u, 6, 9)$, ${\rm lk}(5) = C_7(v, 4, 0, 6$, $9, 1, 8)$, ${\rm
lk}(u) = C_7(8, 6, 4, 3, 7, v, 2)$, ${\rm lk}(6) = C_7(8, u$, $4,
9, 5, 0, 7)$, ${\rm lk}(7) = C_7(8, 6, 0, 1, 3, u, v)$, ${\rm
lk}(8) = C_7(v, 7, 6, u, 2, 1, 5)$ and ${\rm lk}(v) = C_7(7, 8,
5, 4$, $9, 2, u)$. Here $M \cong N_{5}$, by the map $\psi \circ
(0, 2)(3, v, 8, 9, 4)(7, u)$. 

\medskip

\noindent {\bf Case 3.} $(u_1, u_2, u_3) = (3, 9, 6)$, i.e.,
${\rm lk}(1) = C_7(8, 2, 0, 7, 3, 9, 6)$. Now, $02$, $04$, $17$
and $19$ are edges in ${\rm lk}(3)$. Since $237$ and $349$ can
not be faces (because of orientability), ${\rm lk}(3) = C_7(4, 0,
2, 9, 1, 7, x)$ or $C_7(2, 0, 4, 7, 1, 9, y)$, for some $x, y \in
V$. It is not difficult to see that ${\rm lk}(3) = C_7(4, 0, 2,
9, 1, 7, 8)$, $C_7(4, 0, 2, 9, 1, 7, u)$, $C_7(2, 0, 4, 7, 1, 9,
5)$ or $C_7(2, 0, 4, 7, 1, 9, u)$.

If ${\rm lk}(3) = C_7(4, 0, 2, 9, 1, 7, 8)$ then, considering
links of $8$ and $6$ we get $0, 1, 3, 6, 8 \not \in {\rm lk}(v)$
and hence $\deg(v) \leq 6$. So, ${\rm lk}(3) \neq C_7(4, 0, 2, 9,
1, 7, 8)$. If ${\rm lk}(3) = C_7(2, 0, 4, 7, 1, 9, 5)$ then $04,
06, 23$ and $39$ are edges in ${\rm lk}(5)$. Thus $459$ or $256$
is a face. In both the cases, $\deg(v) \leq 6$. Thus, ${\rm lk}(3)
\neq C_7(2, 0, 4, 7, 1, 9, 5)$. By similar arguments, we get
${\rm lk}(3) \neq C_7(4, 0, 2, 9, 1, 7, u)$. So, ${\rm lk}(3) =
(2, 0, 4, 7, 1, 9, u)$.

Now, ${\rm lk}(7) = C_7(4, 3, 1, 0, 6, z, w)$, for some $z, w \in
V$. It is easy to see that $(z, w) = (2, 8)$, $(u, v)$, $(v, 8)$,
$(9, 8)$, $(v, u)$, $(9, v)$.

If $(z, w)= (2, 8)$ then, considering ${\rm lk}(6)$, we get $0, 1,
3, 6, 7 \not \in {\rm lk}(v)$. So, $(z, w) \neq (2, 8)$. If $(z,
w) = (u, v)$ then, considering links of $6, u, 2$ and $8$,
successively, we get $C_3(2, 8, v) \subseteq {\rm lk}(5)$. So,
$(z, w) \neq (u, v)$. If $(z, w) = (v, 8)$ then, $0, 1, 7, 6, 8
\not \in {\rm lk}(u)$. So, $(z, w) \neq (v, 8)$. If $(z, w) = (9,
8)$ then, considering the links of $6, 9$ and $8$ we get $0, 1,
7, 6, 8 \not \in {\rm lk}(u)$. So, $(z, w) \neq (9, 8)$.
Similarly, $(z, w) \neq (v, u)$. Thus, $(z, w) = (9, v)$, i.e.,
${\rm lk}(7) = C_7(4, 3, 1, 0, 6, 9, v)$.

Now, ${\rm lk}(6) = C_7(5, 0, 7, 9, 1, 8, s)$, where, it is easy
to see that $s = u$ or $v$.

\smallskip

\noindent {\bf Subcase 3.1.} $s = u$. Then, completing
successively, we get ${\rm lk}(9) = C_7(3, 1, 6, 7, v, 5, u)$,
${\rm lk}(5) = C_7(0, 4, 8, v, 9, u, 6)$, ${\rm lk}(8) = C_7(1,
2, 4, 5, v, u, 6)$, ${\rm lk}(4) = C_7(0, 3, 7, v, 2, 8, 5)$ and
${\rm lk}(2) = C_7(0, 1, 8, 4, v, u, 3)$. Here, $M \cong N_{6}$
by the map $\psi \circ (0, 9, 3, u)(1, 5, v, 7, 6, 4, 2)$.

\smallskip

\noindent {\bf Subcase 3.2.} $s = v$. Then, completing
successively, we get ${\rm lk}(9) = C_7(1, 3, u, 8, v, 7, 6)$,
${\rm lk}(8) = C_7(1, 2, 5, u, 9, v, 6)$, ${\rm lk}(v) = C_7(4,
7, 9, 8, 6, 5, u)$, ${\rm lk}(u) = C_7(2, 3$, $9, 8, 5, v, 4)$
and ${\rm lk}(2) = C_7(0, 1, 8, 5, 4, u, 3)$. Here, $M \cong
N_{2}$ by the map $\psi \circ (0, 1, 6, 5, 9, u, 3, 7, v, 4, 8)$.

\medskip

\noindent {\bf Case 4.} ${\rm lk}(1) = C_7(8, 2, 0, 7, 3, 9, u)$.
By similar arguments as in the previous cases, ${\rm lk}(3) =
C_7(4, 0, 2, 9, 1, 7, v)$ and ${\rm lk}(2) = C_7(9, 3, 0, 1, 8,
6, 5)$ or $C_7(9, 3, 0, 1, 8, v, 5)$. In the first case, $M \cong
N_{6}$. In the second case, $M \cong N_{4}$.

\medskip

\noindent {\bf Case 5.} ${\rm lk}(1) = C_7(8, 2, 0, 7, 4, 3, 9)$.
By similar arguments as in the first three cases, ${\rm lk}(4) =
C_7(7, 1, 3, 0, 5, 9, u)$, $C_7(7, 1, 3, 0, 5, 8, u)$, $C_7(7, 1,
3, 0, 5, u, 2)$ or $C_7(7, 1, 3, 0, 5, u, v)$.

If ${\rm lk}(4) = C_7(7, 1, 3, 0, 5, 9, u)$ then $M \cong N_{6}$.

If ${\rm lk}(4) = C_7(7, 1, 3, 0, 5, 8, u)$ then ${\rm lk}(3) =
C_7(9, 1, 4, 0, 2, v, 6)$ or $C_7(9, 1, 4, 0, 2, v, u)$. In the
first case, $M \cong N_{5}$. In the second case, $M \cong N_{4}$.

If ${\rm lk}(4) = C_7(7, 1, 3, 0, 5, u, 2)$ then ${\rm lk}(2) =
C_7(8, 1, 0, 3, 7, 4, u)$, ${\rm lk}(3) = C_7(9, 1, 4, 0, 2, 7$,
$v)$, ${\rm lk}(7) = C_7(3, 2, 4, 1, 0, 6, v)$ and ${\rm lk}(6) =
C_7(5, 0, 7, v, u, 9, 8)$ or $C_7(5, 0, 7, v, 8, u, 9)$. In the
first case, $M \cong N_{3}$. In the second case, $M \cong N_{1}$.

If ${\rm lk}(4) = C_7(7, 1, 3, 0, 5, u, v)$ then ${\rm lk}(7) =
C_7(6, 0, 1, 4, v, 8, 2)$ or $C_7(6, 0, 1, 4, v, 9, 8)$. In the
first case, $M \cong N_{6}$. In the second case, $M \cong N_{4}$.

\medskip

\noindent {\bf Case 6.} ${\rm lk}(1) = C_7(8, 2, 0, 7, 4, 6, 5)$.
By similar arguments as in the first three cases, ${\rm lk}(4) =
C_7(6, 1, 7, 5, 0, 3, 9)$. This case is isomorphic to Case 1 by
the map $(4, 1, 0, 6, 2, 5, 3, 9, 8)$.

\medskip

\noindent {\bf Case 7.} ${\rm lk}(1) = C_7(8, 2, 0, 7, 4, 9, 6)$.
By similar arguments as in the first three cases, ${\rm lk}(4) =
C_7(5, 0, 3, 9, 1, 7, u)$ or $C_7(3, 0, 5, 7, 1, 9, u)$. In the
first case, $M \cong N_{4}$. In the second case, $M \cong N_{6}$.

\medskip

\noindent {\bf Case 8.} ${\rm lk}(1) = C_7(8, 2, 0, 7, 4, 9, u)$.
By similar arguments as in the first three cases, ${\rm lk}(4) =
C_7(3, 0, 5, 7, 1, 9, v)$ or $C_7(5, 0, 3, 9, 1, 7, 8)$.

If ${\rm lk}(4) = C_7(3, 0, 5, 7, 1, 9, v)$ then ${\rm lk}(7) =
C_7(6, 0, 1, 4, 5, u, 8)$ or $C_7(6, 0, 1, 4, 5, u, v)$.  In the
first case, $M \cong N_{4}$. In the second case,  $M \cong N_4$
or $N_{5}$.

If ${\rm lk}(4) = C_7(5, 0, 3, 9, 1, 7, 8)$ then ${\rm lk}(8) =
C_7(u, 1, 2, 7, 4, 5, v)$ or $C_7(2, 1, u, 5, 4, 7, v)$. In the
first case, $M \cong N_{6}$. In the second case, $M \cong N_{4}$.

\medskip

\noindent {\bf Case 9.} ${\rm lk}(1) = C_7(8, 2, 0, 7, 5, 4, 9)$.
By similar arguments as in the first three cases, ${\rm lk}(5) =
C_7(7, 1, 4, 0, 6, 8, u)$, $C_7(7, 1, 4, 0, 6, u, 8)$ or $C_7(7,
1, 4, 0, 6, u, v)$.

If ${\rm lk}(5) = C_7(7, 1, 4, 0, 6, 8, u)$ then, ${\rm lk}(8) =
C_7(9, 1, 2, u, 5, 6, 3)$ or $C_7(u, 5, 6, 9, 1, 2, v)$. In the
first case, $M \cong N_{4}$. In the second case, $M \cong N_{5}$.

If ${\rm lk}(5) = C_7(7, 1, 4, 0, 6, u, 8)$ then, ${\rm lk}(8) =
C_7(9, 1, 2, 7, 5, u, 3)$ or $C_7(9, 1, 2, 7, 5, u, v)$. In the
first case, $M \cong N_{6}$. In the second case, $M \cong N_{2}$.

If ${\rm lk}(5) = C_7(7, 1, 4, 0, 6, u, v)$ then ${\rm lk}(7) =
C_7(6, 0, 1, 5, v, 3, 8)$ or $C_7(6, 0, 1, 5, v, 9, 8)$. In the
first case, $M \cong N_{4}$. In the second case, $M \cong N_{6}$.

\medskip

\noindent {\bf Case 10.} ${\rm lk}(1) = C_7(8, 2, 0, 7, 5, 9, u)$.
By similar arguments as in the first three cases, ${\rm lk}(5) =
C_7(7, 1, 9, 4, 0, 6, 8)$, $C_7(7, 1, 9, 4, 0, 6, u)$ or $C_7(7,
1, 9, 4, 0, 6, v)$. In the first case, $M \cong N_{6}$. In the
second case, $M \cong N_{2}$. In the third case, $M \cong N_4$ or
$N_{5}$.

\medskip

\noindent {\bf Case 11.} ${\rm lk}(1) = C_7(8, 2, 0, 7, 9, 4, u)$.
Then, up to an isomorphism, ${\rm lk}(4) = C_7(u, 1, 9, 5, 0$,
$3, 6)$ and ${\rm lk}(6) = C_7(5, 0, 7, 3, 4, u, v)$ or $C_7(7, 0,
5, u, 4, 3, 8)$. In the first case, $M \cong N_{5}$. In the
second case, $M \cong N_{4}$.

\medskip

\noindent {\bf Case 12.} ${\rm lk}(1) = C_7(8, 2, 0, 7, 9, 5,
4)$. By similar arguments as in the first three cases, ${\rm
lk}(5) = C_7(9, 1, 4, 0, 6, 8, u)$ and ${\rm lk}(8) = C_7(2, 1, 4,
6, 5, u, v)$ or $C_7(4, 1, 2, u, 5, 6, v)$. In the first case, $M
\cong N_{6}$. In the second case, $M \cong N_{4}$.

\medskip

\noindent {\bf Case 13.} ${\rm lk}(1) = C_7(8, 2, 0, 7, 9, 5, u)$.
By similar arguments as in the first three cases, ${\rm lk}(5) =
C_7(u, 1, 9, 6, 0, 4, v)$, $C_7(9, 1, u, 4, 0, 6, 3)$ or $C_7(9,
1, u, 4, 0, 6, 8)$. In the first two cases, $M \cong N_{4}$. In
the third case, $M \cong N_{5}$. This completes the proof. \hfill
$\Box$

\bigskip

\noindent {\bf  Acknowledgement\,:} The authors thank the
anonymous referee for many useful references and comments which
helped to improve the presentation of this paper. The second
author thanks CSIR, New Delhi, India for its research fellowship
(Award No.: 9/79(797)/2001\,-\,EMR\,-\,I).

\medskip

\noindent Basudeb Datta and Ashish Kumar Upadhyay, Department of
Mathematics, Indian Institute of Science, Bangalore 560\,012,
India. E-mails: \{dattab, upadhyay\}@math.iisc.ernet.in.

\end{document}